\documentclass[11pt]{article}
\usepackage{graphicx, wrapfig}
\usepackage{flafter}
\usepackage{leftidx}
\usepackage{mathrsfs}
\usepackage{amssymb, amsmath}
\usepackage{bbding}
\usepackage{fancyhdr}
\usepackage{mathrsfs}
\usepackage{color}
\usepackage{stmaryrd}
\usepackage{amsfonts}
\usepackage{latexsym}
\usepackage{psfrag}
\usepackage{graphicx, subfigure}
\usepackage{fancyhdr, graphicx}
\usepackage{multicol}
\usepackage{dsfont}
\usepackage{bbm}
\usepackage{booktabs}
\usepackage[center]{caption2}
\usepackage{cite}
\usepackage{epstopdf}
\usepackage{diagbox}
\usepackage{booktabs}
\usepackage [latin1]{inputenc}
\usepackage{multirow}
\usepackage{enumerate}
\usepackage{enumitem}
\usepackage{algpseudocode, algorithm, algorithmicx}

\allowdisplaybreaks
\renewcommand{\leq}{\leqslant}

\date{}
\newtheorem{theorem}{Theorem}[section]
\newtheorem{lemma}{Lemma}[section]
\newtheorem{remark}{Remark}

\newtheorem{corollary}{Corollary}[section]
\numberwithin{equation}{section}
\newcommand{\zd}{\,\mathrm{d}}

\newcommand{\diff}{\triangledown_{\tau}}
\newcommand{\myvec}[1]{\boldsymbol{#1}}
\newcommand{\abs}[1]{\left|#1\right|}
\newcommand{\absb}[1]{\big|#1\big|}

\newcommand{\abst}[1]{|#1|}
\newcommand{\bra}[1]{\left(#1\right)}
\newcommand{\brab}[1]{\big(#1\big)}
\newcommand{\braB}[1]{\Big(#1\Big)}
\newcommand{\brat}[1]{(#1)}
\newcommand{\kbra}[1]{\left[#1\right]}
\newcommand{\kbrab}[1]{\big[#1\big]}

\newcommand{\myinner}[1]{\left\langle#1\right\rangle}
\newcommand{\myinnerb}[1]{\big\langle#1\big\rangle}
\newcommand{\myinnerB}[1]{\Big\langle#1\Big\rangle}
\newcommand{\mynorm}[1]{\left\|#1\right\|}
\newcommand{\mynormb}[1]{\big\|#1\big\|}

\newcommand{\mynormt}[1]{\|#1\|}
\newcommand{\timenorm}[1]{\absb{\!\absb{\!\absb{#1}\!}\!}}

\begin{document}
\title{Analysis of the second order BDF scheme with variable steps for the molecular beam epitaxial model
without slope selection}
\author{
Hong-lin Liao\thanks{ORCID 0000-0003-0777-6832;
Department of Mathematics, Nanjing University of Aeronautics and Astronautics,
Nanjing 211106, P. R. China. Hong-lin Liao (liaohl@nuaa.edu.cn and liaohl@csrc.ac.cn).
Research supported by a grant 1008-56SYAH18037
from NUAA Scientific Research Starting Fund of Introduced Talent.}
\quad Xuehua Song\thanks{Department of Mathematics, Nanjing University of Aeronautics and Astronautics,
211101, P. R. China.}
\quad Tao Tang\thanks{Department of Mathematics and International Center for Mathematics, Southern
    University of Science and Technology, Shenzhen, Guangdong Province; and
    Division of Science and Technology, BNU-HKBU United International College,
    Zhuhai, Guangdong Province, China.
    Email: tangt@sustech.edu.cn. This author's work is partially supported by the NSF of China under grant number 11731006.}
\quad Tao Zhou\thanks{NCMIS \& LSEC, Institute of Computational Mathematics and Scientific/Engineering Computing,
Academy of Mathematics and Systems Science, Chinese Academy of Sciences, Beijing, 100190,
P. R. China. Email: tzhou@lsec.cc.ac.cn. This author's work is partially supported by the NSF of China (under grant numbers 11822111, 11688101, and 11731006), and the science challenge project (No. TZ2018001).}
}
\maketitle
\normalsize

\begin{abstract}
In this work, we are concerned with the stability and convergence analysis of the second order BDF (BDF2) scheme with variable steps for
the molecular beam epitaxial model without slope selection. We first show that the variable-step BDF2 scheme is convex and uniquely solvable
under a weak time-step constraint. Then we show that it preserves an energy dissipation law if the adjacent time-step ratios $r_k:=\tau_k/\tau_{k-1}<3.561.$ Moreover,
with a novel discrete orthogonal convolution kernels argument and some new estimates on the corresponding positive definite quadratic forms,
the $L^2$ norm stability and rigorous error estimates are established, under the same step-ratios constraint that ensuring the energy stability., i.e., $0<r_k<3.561.$  This is known to be the best result in literature. We finally adopt an adaptive time-stepping strategy to accelerate the computations of  the steady state solution and confirm our theoretical findings by numerical examples.
\end{abstract}


\section{Introduction}
\setcounter{equation}{0}
We consider the following molecular beam epitaxial (MBE) model without slope selection
on a bounded domain $\Omega\subset\mathbb{R}^2$
\begin{align}\label{cont: MBE no-slope problem}
\Phi_t=-\varepsilon\Delta^2\Phi-\nabla\cdot \myvec{f}(\nabla\Phi)
\quad\text{for $\myvec{x}\in\Omega$ and $ 0< t\le T$},
\end{align}
subjected to the initial data $\Phi(\myvec{x}, 0):=\Phi_{0}(\myvec{x})$,  where the nonlinear force vector $\myvec{f}(\myvec{v}):=\frac{\myvec{v}}{1+|\myvec{v}|^2}$. $\Phi=\Phi(\myvec{x}, t)$,
subjected to periodic boundary conditions, is the scaled height function of a thin film in a co-moving frame and
$\varepsilon>0$ is a constant that represents the width of the rounded corners on the otherwise faceted crystalline thin films.

 The above epitaxial growth model admits variable applications in different fields, such as physics \cite{AmarFamily:1996},
biology \cite{EvansThiel:2010} and chemistry \cite{RostKrug:1997}, to name a few.  The MBE model \eqref{cont: MBE no-slope problem},
in which the nonlinear second order term models the Ehrlich-Schwoebel effect and the
linear fourth order term describes the surface diffusion,
defines a gradient flow with respect to the $L^2(\Omega)$ inner product of the following free energy \cite{Golubovic:1997, LiLiu:2003}:
\begin{align}\label{cont: free energy}
E[\Phi]=\int_{\Omega}\kbra{\frac{\varepsilon}{2}(\Delta\Phi)^2
-\frac{1}{2}\ln\brab{1+\abs{\nabla\Phi}^2}}\zd\myvec{x}.
\end{align}
The logarithmic term therein is bounded above by zero
but unbounded below (and has no relative minima),
which implies that no energetically favored values exist for $\nabla\Phi$.
From the physical point of view,
this means that there is no slope selection mechanism.
Thus, it may result in multi-scale behavior in a rough-smooth-rough pattern,
especially at an early stage of epitaxial growth on rough surfaces.
The well-posedness of the initial-boundary-value problem \eqref{cont: MBE no-slope problem}
was studied by Li and Liu in \cite{LiLiu:2003} using the perturbation analysis.
The authors \cite[Theorem 3.3]{LiLiu:2003} proved that,
if the initial data $\phi_0\in H_{per}^{m}(\Omega)$ for some integer $m\ge2$,
the problem has a unique weak solution $\phi$ such that
$\phi\in L^{\infty}\bra{0,T;H^m(\Omega)}\cap L^2\bra{0,T;H^{m+2}(\Omega)}$
and $\partial_t\phi\in L^2\bra{0,T;H^{m-2}(\Omega)}$.
As is well-known, the MBE system \eqref{cont: MBE no-slope problem}
is also volume-conservative, i.e., $\bra{\Phi(t),1}=\bra{\Phi_0,1}$ for $t>0$,
and admits the following energy dissipation law
\begin{align}\label{cont: energy dissipation}
\frac{\zd}{\zd t}E[\Phi]=-\mynorm{\Phi_t}_{L^{2}(\Omega)}^2\le0, \quad  0< t\le T,
\end{align}
where $\bra{\cdot,\cdot}$ denote the inner product in $L^{2}(\Omega)$
and $\mynorm{\cdot}_{L^{2}(\Omega)}$ is the associated norm.
Also, by the Green's formula and Cauchy-Schwarz inequality, one has the $L^2$ norm solution estimate
\begin{align}\label{cont: L2 norm solution estimate}
\mynorm{\Phi}_{L^{2}(\Omega)}\le e^{t/(4\varepsilon)}\mynorm{\Phi_0}_{L^{2}(\Omega)}, \quad  0< t\le T.
\end{align}

As analytic solutions are not in general available, numerical schemes for the above MBE model have been widely studied in recent years
\cite{ChenCondeWangWangWise:2012,ChenWangWang:2014,JuLiQiaoZhang:2018,
QiaoSunZhang:2015,QiaoZhangTang:2011,ShenWangWangWise:2012,XuTang:2006}. Thus include the stabilized semi-implicit scheme \cite{XuTang:2006},
Crank-Nicolson type schemes \cite{QiaoZhangTang:2011}, convex splitting schemes \cite{ChenCondeWangWangWise:2012,ShenWangWangWise:2012}, the exponential time differencing scheme \cite{JuLiQiaoZhang:2018}, to name just a few. The main focus of the above mentioned works were the discrete energy stability, i.e., one constructs a numerical scheme that can inherit the energy dissipation law in the discrete levels.

It is noticed that in all the above mentioned literature, the numerical analysis was performed for uniform time-steps. In this work, we aim at
investigating a nonuniform version of a classic numerical scheme, i.e., the second order BDF (BDF2) scheme with variable time-steps. To this end, we consider the nonuniform time girds
$$0=t_0 \le t_1 \le...\le t_N = T$$
with the time-step sizes $\tau_{k}:=t_{k}-t_{k-1}.$ We denote the maximal step size as $ \tau := \max_{1\le k\le N}\tau_k$ and define the local time-step ratio as $ r_k := \tau_k/\tau_{k-1}$ for $k\ge2$. Given a grid function $v^n=v(t_n)$, we set $ \diff v^n:=v^{n}-v^{n-1}$ and $\partial_{\tau}v^n:=\diff v^n/\tau_n$ for $k\geq1$.
The motivation for using a nonuniform grid is that one can possibly capture the multi-scale beehives in the time domain. However, the numerical analysis for BDF2 with nonuniform grids seems to be highly nontrivial (compared to the uniform-grid case). One few results can be found in literature. For the linear diffusion case, the existing $L^2$ norm stability and error estimates can be found in \cite{Becker:1998,ChenWangYanZhang:2019,Emmrich:2005,LeRoux:1982}. However, for the analysis therein,
the time-step ratio constraint that guarantees the $L^2$ norm stability are always severer
than the classical zero-stability condition $r_k<1+\sqrt{2}$ for ODE problems \cite{CrouzeixLisbona:1984,Grigorieff:1983}.
Moreover, some undesired factors such as $\exp(C_r\Gamma_n)$ appears in the estimates,
where $\Gamma_n$ can be unbounded as the time steps vanish and $C_r$ grows to infinity once the step-ratios approach the zero-stability limit $1+\sqrt{2}$. Associate analysis for nonlinear problems such as the CH equations can be found in \cite{ChenWangYanZhang:2019}. Again, the error estimates therein are presented under the time-step constraints $r_k<1.53$ (worse than the classical zero-stability condition).
 As an exception, in our previous work \cite{LiaoTangZhou:2019maximumAC}, we have presented a novel analysis for the nonuniform BDF2 scheme of the Allen-Cahn equation under the same condition $r_k<1+\sqrt{2}.$ In a very recent work \cite{LiaoZhang:2019linear}, for the linear diffusion problem, the $L^2$ norm stability and convergence estimates are presented under a much improved stability condition
\begin{equation*}
 0< r_k < r_s:=\brat{3+\sqrt{17}}/2 \approx 3.561,  \quad 2\le k \le N.
\end{equation*}
In particular, a novel discrete orthogonal convolution (DOC) kernels argument related to the nonuniform BDF2 scheme is proposed to perform the analysis in \cite{LiaoZhang:2019linear}.
In the current work, we shall pursuit this study for the nonlinear MBE model under the new zero-stability condition.

\subsection{The variable-step BDF2 scheme}

The well known nonuniform BDF2 formula can be expressed as the following convolutional summation
\begin{align}\label{def: BDF2 formula}
D_{2}v^n=\sum_{k=1}^{n}b_{n-k}^{(n)}\diff v^k, \quad  n\ge 1,
\end{align}
 where the discrete convolution kernels $b_{n-k}^{(n)}$ are defined by $b_0^{(1)}:=1/\tau_1$ for $n=1$,
 and for $n \geq 2$ one has
\begin{align}\label{def: BDF2 kernels}
 b_{0}^{(n)}:=\frac{1+2r_n}{\tau_{n}(1+r_n)},
 \quad b_{1}^{(n)}:=-\frac{r_{n}^2}{\tau_{n}(1+r_n)}\quad \text{and}\quad  b_{j}^{(n)}:=0\quad\text{for $j\ge2$}.
 \end{align}
Without loss of generality, we can include the BDF1 formula in \eqref{def: BDF2 formula} by putting $r_1\equiv0$,
and use it to compute the first-level solution for initialization.

To present the fully discrete scheme, for the physical domain $\Omega=(0, L)^2$,
we use a uniform grid with grid lengths $h_x=h_y=h:=L/M$ (with $M$ being an integer)
to yield the discrete domains
$$\Omega_h:=\{\myvec{x}_h=(ih, jh)\,|\,1\le i, j\le M\}, \quad
\textmd{and} \quad  \bar{\Omega}_h:=\{x_h=(ih, jh)\,|\,0\le i, j\le M\}.$$
For the function
$w_h=w(\myvec{x}_h)$, let $$\Delta_{x}w_{ij}:=(w_{i+1, j}-w_{i-1, j})/(2h),
\quad \textmd{and} \quad \delta_{x}^{2}w_{ij}=(w_{i+1, j}-2w_{ij}+w_{i-1, j})/h^2.$$
The operators $\Delta_{y}w_{ij}$ and $\delta_{y}^{2}w_{ij}$ can be defined similarly.
Moreover, the discrete gradient vector and the discrete Laplacian can also be defined accordingly:
$$\nabla_{h}w_{ij}:=(\Delta_{x}w_{ij}, \Delta_{y}w_{ij})^{T},
\quad  \Delta_{h}w_{ij}:=(\delta_x^2+\delta_y^2)w_{ij}.$$
One can further define the discrete divergence as
$\nabla_{h}\cdot \myvec{u}_{ij}:=\Delta_{x}v_{ij}+\Delta_{y}w_{ij}$
for the vector $\myvec{u}_h=(v_h, w_h)^T$.
We also denote the space of $L$-periodic grid functions as
$$\mathbb{V}_h:=\big\{v_h\,|\,v_h\;\text{is $L$-periodic for}\;\myvec{x}_h\in\bar{\Omega}_h\big\}.$$

We are now ready to present the fully implicit variable-step BDF2 scheme
for the MBE equation \eqref{cont: MBE no-slope problem}:
find the numerical solution $\phi_{h}^n\in \mathbb{V}_h$ such that $\phi_{h}^0:=\Phi_{0}(\myvec{x}_{h})$ and
\begin{align}\label{scheme: BDF2 implicit MBE-no-slope}
D_{2}\phi_{h}^n+\varepsilon\Delta_{h}^2\phi_{h}^n+\nabla_{h}\cdot \myvec{f}(\nabla_{h}\phi_{h}^n)=0
\quad\text{for $\myvec{x}_{h}\in\Omega_h$ and $1\le n\le N$.}
\end{align}

\subsection{Summary of the main contributions}
As mentioned, we shall pursuit the further study of the analysis technique in \cite{LiaoZhang:2019linear} for the nonlinear MBE model. In \cite{LiaoZhang:2019linear}, the discrete orthogonal convolution (DOC) kernels are proposed for analyzing the linear diffusion problems. The DOC kernels are defined as follows
\begin{align}\label{def: DOC-Kernels}
\theta_{0}^{(n)}:=\frac{1}{b_{0}^{(n)}}
\quad \mathrm{and} \quad
\theta_{n-k}^{(n)}:=-\frac{1}{b_{0}^{(k)}}
\sum_{j=k+1}^n\theta_{n-j}^{(n)}b_{j-k}^{(j)}
\quad \text{for $1\le k\le n-1$}.
\end{align}
It is easy to verify that the following discrete orthogonal identity holds
\begin{align}\label{eq: orthogonal identity}
 \sum_{j=k}^n\theta_{n-j}^{(n)}b_{j-k}^{(j)}\equiv \delta_{nk}\quad\text{for $1\le k\le n$,}
 \end{align}
where $\delta_{nk}$ is the Kronecker delta symbol.

The main motivation for introducing the DOC kernels lies in the following equality
  \begin{align}\label{eq: orthogonal equality for BDF2}
  \sum_{j=1}^n\theta_{n-j}^{(n)}D_2v^j=\diff v^n\quad\text{for $n\ge1$,}
   \end{align}
which can be derived by exchanging the summation order and using the identity \eqref{eq: orthogonal identity}.

In this work, by showing some new properties of the DOC kernels
and the corresponding quadratic forms (see Lemmas \ref{lem: DOC matrix Quadratic inequality}--\ref{lem: DOC Quadratic inequality-nonlinear}), we are able to show the energy stability and a rigorous error estimate of the nonuniform BDF2 scheme for the nonlinear  MBE model, under the following mild time-step ratios constraint
\begin{enumerate}[itemindent=1em]
\item[\textbf{S1}.]
 $0< r_k < r_s:=\bra{3+\sqrt{17}}/2 \approx 3.561$ for $2\le k \le N$.
\end{enumerate}
This coincides with the results in the linear case \cite{LiaoZhang:2019linear}, and up to now seems to be the best results for nonlinear problems in literature.

The rest of this paper is organized as follows. In the next section, we show that the solution of nonuniform BDF2 scheme is equivalent to the minimization problem
of a convex energy functional, thus it is uniquely solvable. Then, we present in Theorem \ref{thm: energy stable} a discrete energy dissipation law.
In Section 3, we present some new properties of the DOC kernels. This is used in Section 4 to show the $L^2$ norm stability and convergence property of the fully implicit scheme Numerical experiments are presented in Section 5 to show the effectiveness of the BDF2 scheme with an adaptive time-stepping strategy. We finally give some concluding remarks in Section 6.

\section{Solvability and energy stability}
\setcounter{equation}{0}
In this section, we show the solvability and discrete energy stability.
To this end, for any grid functions $v, w\in \mathbb{V}_h$,
we define the discrete $L^2$ inner product $\myinner{v, w}:=h^2\sum_{\myvec{x}_h\in\Omega_h}v_hw_h$
and the associated $L^2$ norm $\mynorm{v}:=\sqrt{\myinner{v, v}}$.
The discrete seminorms $\|\nabla_hv\|$ and $\mynorm{\Delta_hv}$ can be defined respectively by
$$\mynorm{\nabla_hv}:=\sqrt{h^2\sum_{\myvec{x}_h\in\Omega_h}\abs{\nabla_hv_h}^2}
\quad\text{and}\quad
\mynorm{\Delta_hv}:=\sqrt{h^2\sum_{\myvec{x}_h\in\Omega_h}\abs{\Delta_hv_h}^2}\quad\text{for $v\in \mathbb{V}_h$.}$$
For any grid functions $v,w\in\mathbb{V}_h,$ the discrete Green's formula with periodic boundary conditions yield
$\myinner{-\nabla_h\cdot\nabla_hv, w}=\myinner{\nabla_hv, \nabla_hw}$. 
It is easy to verify that for $\epsilon>0$ and $v\in \mathbb{V}_h$
\begin{align}\label{ieq: H1 embeds H2}
\mynorm{\nabla_hv}^2\le\myinner{-\Delta_h v,v}\le \mynorm{\Delta_hv}\cdot\mynorm{v}\le
\frac{\epsilon}{2}\mynorm{\Delta_hv}^2+\frac1{2\epsilon}\mynorm{v}^2.
\end{align}

\subsection{Unique solvability}

We first show the solvability of the BDF2 scheme \eqref{scheme: BDF2 implicit MBE-no-slope}
via a discrete energy functional $G$ on the space $\mathbb{V}_h$,
\begin{align*}
G[z]:=\frac{1}{2}b_0^{(n)}\mynorm{z-\phi^{n-1}}^2+b_1^{(n)}\myinner{\diff\phi^{n-1}, z}+\frac{\varepsilon}{2}\mynorm{\Delta_hz}^2
-\frac{1}{2}\myinnerb{\ln(1+\abs{\nabla_hz}^2), 1}.
\end{align*}
We have the following theorem:
\begin{theorem}\label{thm: uniquely solvable}
If the time-step sizes $\tau_n\le4\varepsilon$, the BDF2 time-stepping scheme \eqref{scheme: BDF2 implicit MBE-no-slope} is
convex \cite{XuLiWuBousquet:2019} and thus uniquely solvable.
\end{theorem}
\begin{proof}
To handle the logarithmic term in the above discrete energy functional $G$, we consider a function $g(\lambda):=\frac{1}{2}\ln\brab{1+\abs{\myvec{u}+\lambda \myvec{v}}^2}$ for any vectors $\myvec{u},\myvec{v}$ such that
\begin{align*}
\frac{\zd g(\lambda)}{\zd\lambda}\Big|_{\lambda=0}=
\frac{\myvec{v}^T\myvec{u}}{1+\abs{\myvec{u}}^2}=\myvec{v}^T\myvec{f}\brat{\myvec{u}}\quad\text{and}\quad
\frac{\zd^2g(\lambda)}{\zd\lambda^2}\Big|_{\lambda=0}=
\frac{1-\abs{\myvec{u}}^2}{\brat{1+\abs{\myvec{u}}^2}^2}\myvec{v}^T\myvec{v}\le\myvec{v}^T\myvec{v}\,.
\end{align*}
For any time-level index $n\geq1$, the time-step condition implies $b_0^{(n)}>\frac{1}{4\varepsilon}$.
Then the functional $G$ is strictly convex as for any $\lambda\in\mathbb{R}$ and any $\psi_h\in\mathbb{V}_h$, one has
\begin{align*}
\frac{\zd^2}{\zd\lambda^2}G[z+\lambda\psi]\Big|_{\lambda=0}
\geq&\,  b_0^{(n)}\mynorm{\psi}^2+\varepsilon\mynorm{\Delta_h\psi}^2
-\mynorm{\nabla_h\psi}^2
\geq (b_0^{(n)}-\frac{1}{4\varepsilon})\mynorm{\psi}^2\geq0,
\end{align*}
where the inequality \eqref{ieq: H1 embeds H2} with $\epsilon:=2\varepsilon$ was applied to bound $\mynorm{\nabla_h\psi}^2$
in the above derivation.
Thus, the functional $G$ admits a unique minimizer (denoted by $\phi_h^n$) if and only if it solves
\begin{align*}
0=&\,\frac{\zd}{\zd\lambda}G[z+\lambda\psi]\Big|_{\lambda=0}\\
=&\,b_0^{(n)}\myinner{z-\phi^{n-1}, \psi}+b_1^{(n)}\myinner{\diff\phi^{n-1}, \psi}+\varepsilon\myinner{\Delta_h^2z, \psi}
-\myinnerb{\myvec{f}\brat{\nabla_hz}, \nabla_h\psi}\\
=&\, \myinnerB{b_0^{(n)}(z-\phi^{n-1})+b_1^{(n)}\diff\phi^{n-1}+\varepsilon\Delta_h^2z+\nabla_h\cdot\myvec{f}\brat{\nabla_hz}, \psi}.
\end{align*}
This equation holds for any $\psi_h\in\mathbb{V}_h$ if and only if the unique minimizer $\phi_h^n\in\mathbb{V}_h$ solves
\begin{align*}
b_0^{(n)}(\phi_h^n-\phi_h^{n-1})+b_1^{(n)}\diff\phi_h^{n-1}+\varepsilon\Delta_h^2\phi_h^n
+\nabla_h\cdot \myvec{f}\bra{\nabla_h\phi_h^n}=0,
\end{align*}
and this coincides with the BDF2 scheme \eqref{scheme: BDF2 implicit MBE-no-slope}. The proof is completed.
\end{proof}

\subsection{Discrete energy dissipation law}

To establish the energy stability of the BDF2 scheme \eqref{scheme: BDF2 implicit MBE-no-slope}, we first present the following lemma for which the proof is similar as in \cite[Lemma 2.1]{LiaoZhang:2019linear}.
\begin{lemma}\label{lem: BDF2 kernel positiveDef}
Suppose that \textbf{S1} holds, then for any non-zero sequence $\{w_k\}_{k=1}^n,$ it holds
\begin{align}\label{Convolution kernel inequality}
2w_{k}\sum_{j=1}^{k}b_{k-j}^{(k)}w_j\ge&\, \frac{r_{k+1}}{1+r_{k+1}}\frac{w_k^2}{\tau_k}-\frac{r_k}{1+r_k}\frac{w_{k-1}^2}{\tau_{k-1}}
+\braB{\frac{2+4r_k-r_{k}^2}{1+r_k}-\frac{r_{k+1}}{1+r_{k+1}}}\frac{w_k^2}{\tau_k}, \,\, k\ge 2.
\end{align}
Consequently, the discrete convolution kernels $b_{n-k}^{(n)}$ are positive definite in the sense that
\begin{align*}
\sum_{k=1}^nw_k\sum_{j=1}^kb_{k-j}^{(k)}w_j\ge \frac1{2}\sum_{k=1}^n\braB{\frac{2+4r_k-r_{k}^2}{1+r_k}-\frac{r_{k+1}}{1+r_{k+1}}}\frac{w_k^2}{\tau_k}>0, \quad n\ge 2.
\end{align*}
\end{lemma}

Notice that the BDF2 formula \eqref{def: BDF2 formula} is a multi-step scheme, thus it is nature to consider the following modified discrete energy
\begin{align*}
\mathcal{E}[\phi^n]:=E[\phi^n]+\frac{r_{n+1}}{2(1+r_{n+1})\tau_n}\mynormb{\diff\phi^n}^2, \quad 0\le n\le N,
\end{align*}
where $\mathcal{E}[\phi^0]=E[\phi^0]$ due to $r_1\equiv0$ and
$E[\phi^n]$ is the discrete version of the energy functional \eqref{cont: free energy}, i.e.,
\begin{align}\label{def: discrete free energy}
E[\phi^n]:=\frac{\varepsilon}{2}\mynormb{\Delta _{h}\phi^n}^2-\frac{1}{2}\myinnerb{\ln(1+\abs{\nabla_h\phi^n}^2), 1}\quad\text{for $0\le n\le N$.}
\end{align}
To establish an energy dissipation law, we impose a restriction of time-step sizes $\tau_n$ as follows
\begin{align}\label{ieq: time-steps restriction}
\tau_n\le4\varepsilon\min\bigg\{1,\frac{2+4r_n-r_n^2}{1+r_n}-\frac{r_{n+1}}{1+r_{n+1}}\bigg\}, \quad n\geq1.
\end{align}
We are now ready to present the following theorem.
\begin{theorem}\label{thm: energy stable}
Suppose that \textbf{S1} holds with the time-step condition \eqref{ieq: time-steps restriction}, then the BDF2 scheme \eqref{scheme: BDF2 implicit MBE-no-slope} admits the following  energy dissipation law:
\begin{align*}
\mathcal{E}[\phi^n]\le \mathcal{E}[\phi^{n-1}]\le\mathcal{E}[\phi^{0}]=E[\phi^{0}], \quad n\geq1.
\end{align*}
\end{theorem}
\begin{proof}
Taking the inner product of \eqref{scheme: BDF2 implicit MBE-no-slope} by $\diff\phi^n$,  one has
\begin{align}\label{thmProof: energy inners}
\myinnerb{D_{2}\phi^n , \diff\phi^n}
+\varepsilon\myinnerb{\Delta_{h}^2\phi^n,  \diff\phi^n}
+\myinnerb{\nabla_{h}\cdot \myvec{f}(\nabla_h\phi^n), \diff\phi^n}=0\quad\text{for $n\geq1$.}
\end{align}
By using the summation by parts argument and $2a(a-b)=a^2-b^2+(a-b)^2$, we obtain
\begin{align}\label{thmProof: energy-inner1}
\varepsilon\myinner{\Delta_{h}^2\phi^n, \diff\phi^n}=\varepsilon\myinner{\Delta_{h}\phi^n, \Delta_{h}\diff\phi^n}
=\frac{\varepsilon}{2}\mynormb{\Delta_{h}\phi^n}^2-\frac{\varepsilon}{2}\mynormb{\Delta_{h}\phi^{n-1}}^2
+\frac{\varepsilon}{2}\mynormb{\Delta_{h}\diff\phi^n}^2.
\end{align}
To deal with the nonlinear term at the left-hand of \eqref{thmProof: energy inners}, we notice that for any vectors $\myvec{u},\myvec{v}$
one has
\begin{align*}
\frac{2(\myvec{u}-\myvec{v})^T\myvec{u}}{1+\abs{\myvec{u}}^2}
=&\,\frac{\abs{\myvec{u}}^2-\abs{\myvec{v}}^2}{1+\abs{\myvec{u}}^2}
+\frac{\abs{\myvec{u}-\myvec{v}}^2}{1+\abs{\myvec{u}}^2}
\le\ln\frac{1+\abs{\myvec{u}}^2}{1+\abs{\myvec{v}}^2}
+\abs{\myvec{u}-\myvec{v}}^2,
\end{align*}
where the inequality $\frac{z}{1+z}\le \ln(1+z)$  with
$z=\brat{\abs{\myvec{u}}^2-\abs{\myvec{v}}^2}/\brat{1+\abs{\myvec{v}}^2}>-1$ was used. Thus,
by taking $\myvec{u}:=\nabla_h\phi^n$ and $\myvec{v}:=\nabla_h\phi^{n-1}$, one has
\begin{align}\label{thmProof: energy-inner2}
\myinnerb{\nabla_{h}\cdot \myvec{f}(\nabla_h\phi^n), \diff\phi^n}
=&\, -\myinnerb{\myvec{f}(\nabla_h\phi^n), \nabla_{h}\diff\phi^n}\nonumber\\
\ge&\, -\frac{1}{2}\myinnerb{\ln(1+\abst{\nabla_h\phi^n}^2), 1}+\frac{1}{2}\myinnerb{\ln(1+\abst{\nabla_h\phi^{n-1}}^2), 1}
-\frac{1}{2}\mynormb{\nabla_h\diff\phi^{n}}^2\nonumber\\
\ge&\, -\frac{1}{2}\myinnerb{\ln(1+\abst{\nabla_h\phi^n}^2), 1}+\frac{1}{2}\myinnerb{\ln(1+\abst{\nabla_h\phi^{n-1}}^2), 1}\nonumber\\
&\, -\frac{\varepsilon}{2}\mynormb{\Delta_h\diff\phi^{n}}^2-\frac{1}{8\varepsilon}\mynormb{\diff\phi^{n}}^2,
\end{align}
where the inequality \eqref{ieq: H1 embeds H2} with $v:=\diff\phi^{n}$ and $\epsilon:=2\varepsilon$ was applied to bound $\mynormb{\nabla_h\diff\phi^{n}}^2$ in the last step. By inserting
\eqref{thmProof: energy-inner1}-\eqref{thmProof: energy-inner2} into \eqref{thmProof: energy inners} and using together
the definition \eqref{def: discrete free energy}, we obtain
\begin{align}\label{thmProof: energy inequality}
\myinnerb{D_{2}\phi^n,\diff\phi^n}-\frac{1}{8\varepsilon}\mynormb{\diff\phi^{n}}^2+E[\phi^n]-E[\phi^{n-1}]\le0\quad\text{for $n\geq1$.}
\end{align}
We now proceed the proof by dealing with the first term at the left-hand of \eqref{thmProof: energy inequality}.
For $n\ge2$, Lemma \ref{lem: BDF2 kernel positiveDef} and the time-step condition \eqref{ieq: time-steps restriction} yield
\begin{align*}
\myinner{D_{2}\phi^n , \diff\phi^n}
\geq&\, \frac{r_{n+1}}{2(1+r_{n+1})\tau_n}\mynormb{\diff\phi^n}^2-
\frac{r_{n}}{2(1+r_{n})\tau_{n-1}}\mynormb{\diff\phi^{n-1}}^2+\frac{1}{8\varepsilon}\mynormb{\diff\phi^n}^2.
\end{align*}
Then it follows from \eqref{thmProof: energy inequality} that
$$\mathcal{E}[\phi^n]\le \mathcal{E}[\phi^{n-1}], \quad n\ge2.$$
For the case $n=1$, the facts $r_1=0$ and the time-step condition \eqref{ieq: time-steps restriction} yield
$\tau_1\le\frac{4\varepsilon(2+r_{2})}{1+r_{2}}.$ Consequently, one has
\begin{align*}
\myinner{D_{2}\phi^1 , \diff\phi^1}=\frac{1}{\tau_1}\mynormb{\diff\phi^1}^2
\geq&\, \frac{r_{2}}{2(1+r_{2})\tau_1}\mynormb{\diff\phi^1}^2+\frac{1}{8\varepsilon}\mynormb{\diff\phi^1}^2.
\end{align*}
By inserting the above inequality into \eqref{thmProof: energy inequality}, one gets
$$\mathcal{E}[\phi^1]\le E[\phi^{0}]=\mathcal{E}[\phi^{0}].$$
This completes the proof.
\end{proof}
\begin{remark}\label{remark: modified discrete energy}
Obviously, $\mathcal{E}[\phi^n]-E[\phi^n]\approx\tau_n\mynormt{\partial_{\tau}\phi^n}^2$ so that the modified energy approximates
the original energy with an order of $O(\tau_n)$. From the computational view of point, the modified discrete energy form
$\mathcal{E}[\phi^n]$ suggests that small time-steps
(with small step ratios) are necessary to capture the solution behaviors when
$\mynormt{\partial_t\phi}$ becomes large, while large time-steps are acceptable to accelerate
the time integration when $\mynormt{\partial_t\phi}$ is small.
\end{remark}
\begin{remark}\label{remark: time-steps restriction}
Notice that the first time-step condition in \eqref{ieq: time-steps restriction} comes from the unique solvability
and the second one  is necessary to maintain the discrete energy stability. In practice, the time-step constraint \eqref{ieq: time-steps restriction} requires $\tau_n=O(\varepsilon)$ which is essentially determined by the value of surface diffusion parameter $\varepsilon$.
Thus, the time-step condition is acceptable since the restriction $\tau_n=O(\varepsilon)$
is always required in the $L^2$ norm stability or convergence analysis \cite{QiaoSunZhang:2015,ChenCondeWangWangWise:2012,ChenWangWang:2014,JuLiQiaoZhang:2018}.
\end{remark}

\section{New properties of the DOC kernels}
\setcounter{equation}{0}

We firstly present some basic properties of the DOC kernels which can be found
in \cite[Lemma 2.2, Corollary 2.1 and Lemma 2.3]{LiaoZhang:2019linear}.
\begin{lemma}\label{lem: DOC property}
	Under the assumption \textbf{S1}, which implies that the discrete convolution kernels  $b^{(n)}_{n-k}$ in
	\eqref{def: BDF2 kernels} are positive semi-definite, then the following properties of the DOC kernels $\theta_{n-j}^{(n)}$ hold:
\begin{itemize}
  \item[(I)] The discrete kernels $\theta_{n-j}^{(n)}$ are positive definite;
  \item[(II)] The discrete kernels $\theta_{n-j}^{(n)}$ are positive and $\displaystyle \theta_{n-j}^{(n)}=\frac{1}{b^{(j)}_{0}}\prod_{i=j+1}^n\frac{r_i^2}{1+2r_i}$ for $1\le j\le n$;
  \item[(III)] $\displaystyle \sum_{j=1}^{n}\theta_{n-j}^{(n)}=\tau_n$ such that
  $\displaystyle \sum_{k=1}^{n}\sum_{j=1}^{k}\theta_{k-j}^{(k)}=t_n$ for $n\ge1$.
\end{itemize}
\end{lemma}

In order to facilitate the following numerical analysis,
we use the BDF2 kernels $b_{k-j}^{(k)}$, the DOC kernels $\theta_{k-j}^{(k)}$ and the $2\times 2$ identity matrix $\mathbf{I}_2$
to define the following matrices
\begin{equation*}
\mathbf{B}_2:=\left(
\begin{array}{cccc}
b_0^{(1)} &  & &\\
b_1^{(2)} & b_0^{(2)} & &\\
 & \ddots &\ddots & \\
  &  &b_1^{(n)} &b_0^{(n)} \\
\end{array}
\right)\otimes\mathbf{I}_2, \qquad
\mathbf{\Theta}_2:=\left(
\begin{array}{cccc}
\theta_0^{(1)} &  & &\\
\theta_1^{(2)} & \theta_0^{(2)} & &\\
\vdots& \vdots &\ddots & \\
\theta_{n-1}^{(n)} &\theta_{n-2}^{(n)} &\ldots &\theta_0^{(n)} \\
\end{array}
\right)\otimes\mathbf{I}_2,
\end{equation*}
where ``$\otimes$" denotes the tensor product.
By the discrete orthogonal identity \eqref{eq: orthogonal identity},
one can verify that $\mathbf{\Theta}_2 = \mathbf{B}_2^{-1}.$
Lemma \ref{lem: BDF2 kernel positiveDef} show that the real symmetric matrix
\begin{align}\label{def:matrix B}
\mathbf{B}: = \mathbf{B}_2+\mathbf{B}_2^T\quad\text{is positive definite.}
\end{align}
Similarly, Lemma \ref{lem: DOC property}(I) implies that the real symmetric matrix
$\mathbf{\Theta}: = \mathbf{\Theta}_2+\mathbf{\Theta}_2^T$ is positive definite.
By using \eqref{def:matrix B}, one can check that
\begin{align}\label{def:matrix Theta formula}
\mathbf{\Theta}=\mathbf{B}_2^{-1}+(\mathbf{B}_2^{-1})^T=(\mathbf{B}_2^{-1})^T\mathbf{B}\mathbf{B}_2^{-1}.
\end{align}
Moreover, we define a diagonal matrix
$\mathbf{\Lambda}_{\tau}:=\mathrm{diag}(\sqrt{\tau_1},\sqrt{\tau_2},\ldots,\sqrt{\tau_n})\otimes\mathbf{I}_2$ and
\begin{align}\label{def:matrix tilde B}
\widetilde{\mathbf{B}}_2:=\mathbf{\Lambda}_{\tau}\mathbf{B_2}\mathbf{\Lambda}_{\tau}=\left(
\begin{array}{cccc}
\tilde{b}_0^{(1)} &  & &\\
\tilde{b}_1^{(2)}& \tilde{b}_0^{(2)} & &\\
 & \ddots &\ddots & \\
  &  &\tilde{b}_1^{(n)} &\tilde{b}_0^{(n)}\\
\end{array}
\right)\otimes\mathbf{I}_2
\end{align}
where the discrete kernels $\tilde{b}_0^{(k)}$ and
$\tilde{b}_1^{(k)}$ are given by $(r_1\equiv0)$
\begin{align*}
\tilde{b}_0^{(k)}=\frac{1+2r_k}{1+r_k}\quad \mathrm{and} \quad \tilde{b}_1^{(k)}=-\frac{{r_k}^{3/2}}{1+r_k}\qquad \text{for $1\le k\le n$.}
\end{align*}
By following the proof of \cite[Lemma A.1]{LiaoJiZhang:2020},
it is easy to check that the real symmetric matrix
\begin{align*}
\widetilde{\mathbf{B}}:=\widetilde{\mathbf{B}}_2+\widetilde{\mathbf{B}}_2^T
\quad\text{is positive definite.}
\end{align*}
So there exists a non-singular upper triangular matrix $\mathbf{L}$ such that
\begin{align}\label{decom:mateix B}
\widetilde{\mathbf{B}}=\mathbf{\Lambda}_{\tau}
\mathbf{B}\mathbf{\Lambda}_{\tau}
=\mathbf{L}^T\mathbf{L}\quad \text{or}\quad
\mathbf{B}=(\mathbf{L}\mathbf{\Lambda}_{\tau}^{-1})^T
\mathbf{L}\mathbf{\Lambda}_{\tau}^{-1}.
\end{align}

We will present some discrete convolution inequalities with respect to the DOC kernels.
To do so, we introduce the vector norm $\timenorm{\cdot}$ by
$\timenorm{\myvec{u}}:=\sqrt{\myvec{u}^T\myvec{u}}$ and
the associated matrix norm
$\timenorm{\textbf{A}}:=\sqrt{\rho\brab{\textbf{A}^T\textbf{A}}}$.
Also, define a positive quantity
\begin{align}\label{def: quantity Mrn}
\mathcal{M}_r:=\max_{n\ge1}\timenorm{\widetilde{\mathbf{B}}_2}^2\timenorm{\mathbf{L}^{-1}}^4
=\max_{n\ge1}\frac{\lambda_{\max}\brab{\widetilde{\mathbf{B}}_2^T\widetilde{\mathbf{B}}_2}}{\lambda_{\min}^2\brab{\widetilde{\mathbf{B}}}}.
\end{align}
Under the step-ratio condition \textbf{S1},
a rough estimate $\mathcal{M}_r<39$ could be followed from \cite[Lemmas A.1 and A.2]{LiaoJiZhang:2020}.
As noticed in \cite[Remark 3]{LiaoJiZhang:2020}, one has $\mathcal{M}_r\le4$ if
practical simulations do not continuously use large step-ratios approaching the stability limit $r_s=3.561$.

\begin{lemma}\label{lem: DOC matrix Quadratic inequality}
If the condition \textbf{S1} holds, then for any vector sequences
$\myvec{z}^k, \myvec{w}^k\in \mathbb{R}^2$ $(1\le k\le n)$,
\begin{align*}
\sum_{k=1}^n\sum_{j=1}^k\theta_{k-j}^{(k)}\brab{\myvec{z}^k}^{T}\myvec{w}^j
\le&\,\frac{\epsilon}{2} \myvec{z}^{T}\mathbf{\Theta} \myvec{z}+\frac{1}{2\epsilon}\myvec{w}^{T}\mathbf{B}^{-1}\myvec{w}
\quad \text{for any $\epsilon\ge0$.}
\end{align*}
where the vector $\myvec{z}:=\brab{\brat{\myvec{z}^1}^T,\brat{\myvec{z}^2}^T,\cdots,\brat{\myvec{z}^n}^T}^T$
and $\myvec{w}:=\brab{\brat{\myvec{w}^1}^T,\brat{\myvec{w}^2}^T,\cdots,\brat{\myvec{w}^n}^T}^T$.
\end{lemma}
\begin{proof}
This result can be verified by following from the proof of \cite[Lemma A.3]{LiaoJiZhang:2020}.
\end{proof}


\begin{lemma}\cite[Lemma 3.5]{JuLiQiaoZhang:2018}\label{lem: lipschitz property-nonlinear}
For any $\myvec{v}, \myvec{w}\in \mathbb{R}^2$, there exists a symmetric matrix $\mathbf{Q}_f\in \mathbb{R}^{2\times2}$ such that
$\myvec{f}(\myvec{v})-\myvec{f}(\myvec{w})=\mathbf{Q}_f(\myvec{v}-\myvec{w})$,
and the eigenvalues of $\mathbf{Q}_f$ satisfy $\lambda_1, \lambda_2\in[-1/8,1]$.
Consequently, it holds that
\begin{align*}
\abs{\myvec{f}(\myvec{v})-\myvec{f}(\myvec{w})}\le\abs{\myvec{v}-\myvec{w}}\quad\text{for any $\myvec{v}, \myvec{w}\in \mathbb{R}^2$.}
\end{align*}
\end{lemma}

\begin{lemma}\label{lem: DOC Quadratic inequality-nonlinear}
Assume that the condition \textbf{S1} holds. For any vector sequences $\myvec{v}^k, \myvec{z}^k, \myvec{w}^k\in \mathbb{R}^2$, $1\le k\le n$
 and any $\epsilon>0$, it holds that
\begin{align*}
\sum_{k=1}^{n}\sum_{j=1}^{k}\theta_{k-j}^{(k)}\brat{\myvec{z}^k}^T
\kbra{\myvec{f}(\myvec{v}^j+\myvec{w}^j)-\myvec{f}(\myvec{v}^j)}
\le \sum_{k=1}^{n}\sum_{j=1}^{k}\theta_{k-j}^{(k)}\kbra{\epsilon\brat{\myvec{z}^k}^T\myvec{z}^j
+\frac{\mathcal{M}_r}{\epsilon}\brat{\myvec{w}^k}^T\myvec{w}^j},
\end{align*}
where the positive constant $\mathcal{M}_r$, independent of the time $t_n$,
time-step sizes $\tau_n$ and time-step ratios $r_n$, is defined by \eqref{def: quantity Mrn}. Consequently,
\begin{align*}
\sum_{k=1}^{n}\sum_{j=1}^{k}\theta_{k-j}^{(k)}\brat{\myvec{z}^k}^T
\kbra{\myvec{f}(\myvec{v}^j+\myvec{z}^j)-\myvec{f}(\myvec{v}^j)}
\le 2\sqrt{\mathcal{M}_r}\sum_{k=1}^{n}\sum_{j=1}^{k}\theta_{k-j}^{(k)}\brat{\myvec{z}^k}^T\myvec{z}^j.
\end{align*}
\end{lemma}

\begin{proof}According to Lemma \ref{lem: lipschitz property-nonlinear},
there exists a sequence of symmetric matrices
$$\mathbf{Q}_f^j\in \mathbb{R}^{2\times2} \quad\text{such that}\quad
\myvec{f}(\myvec{w}^j+\myvec{v}^j)-\myvec{f}(\myvec{w}^j)
=\mathbf{Q}_f^j\myvec{v}^j\quad\text{for $1\le j\le n$,}$$
where the corresponding eigenvalues of $\mathbf{Q}_f^j$ satisfy
$\lambda_{j1}, \lambda_{j2}\in[-1/8,1]$ for $1\le j\le n$.
Now we define the following symmetric matrix
$$\mathbf{Q}:=\mathrm{diag}\brat{\mathbf{Q}_f^1,
\mathbf{Q}_f^2,\cdots,\mathbf{Q}_f^n}\in \mathbb{R}^{2n\times2n}.$$
The eigenvalues $\mu_{k}$ of $\mathbf{Q}_{2n\times2n}$ satisfy
$\mu_{k}\in[-1/8,1]$ for $1\le k\le 2n$. Thus
\begin{align}\label{lemproof: Quadratic inequality1}
\rho\brab{\mathbf{Q}}\le1\quad\text{such that}\quad \timenorm{\mathbf{Q}}\le1.
\end{align}
Also, it is easy to verify that $\mathbf{Q}$ and $\mathbf{\Lambda}_{\tau}$ are commutative, that is, $\mathbf{Q}\mathbf{\Lambda}_{\tau}=\mathbf{\Lambda}_{\tau}\mathbf{Q}$.

By introducing $\myvec{z}_e:=\brab{\brat{\myvec{z}^1}^T,\brat{\myvec{z}^2}^T,\cdots,\brat{\myvec{z}^n}^T}^T$
and $\myvec{w}_e:=\brab{\brat{\myvec{w}^1}^T,\brat{\myvec{w}^2}^T,\cdots,\brat{\myvec{w}^n}^T}^T$,
we apply Lemma \ref{lem: DOC matrix Quadratic inequality} with
$\myvec{z}:=\myvec{z}_e$ and $\myvec{w}:=\mathbf{Q}\myvec{w}_e$ to derive that
\begin{align}\label{lemProof: DOC Quadratic inequality-nonlinear}
\sum_{k=1}^{n}\sum_{j=1}^{k}\theta_{k-j}^{(k)}\brat{\myvec{z}^k}^T
&\,\kbra{\myvec{f}(\myvec{v}^j+\myvec{w}^j)-\myvec{f}(\myvec{v}^j)}
=\sum_{k=1}^{n}\sum_{j=1}^{k}\theta_{k-j}^{(k)}
\brat{\myvec{z}^k}^T\mathbf{Q}_f^j\myvec{w}^j\nonumber\\
\le&\,\frac{\epsilon}{2}\myvec{z}_e^{T}\mathbf{\Theta}\myvec{z}_e
+\frac1{2\epsilon}\myvec{w}_e^{T}\mathbf{Q}^T\mathbf{B}^{-1}\mathbf{Q}\myvec{w}_e\nonumber\\
=&\,\epsilon\sum_{k=1}^{n}\sum_{j=1}^{k}
\theta_{k-j}^{(k)}\brat{\myvec{z}^k}^T\myvec{z}^j
+\frac1{2\epsilon}\myvec{w}_e^{T}\mathbf{Q}^T\mathbf{B}^{-1}\mathbf{Q}\myvec{w}_e.
\end{align}
Now we deal with the second term at the right side of the above inequality.
It follows from \eqref{def:matrix Theta formula} and \eqref{decom:mateix B} that
\begin{align*}
\mathbf{\Theta}=(\mathbf{B}_2^{-1})^T\mathbf{B}\mathbf{B}_2^{-1}
=(\mathbf{B}_2^{-1})^T\bra{\mathbf{L}\mathbf{\Lambda}_\tau^{-1}}^T
\mathbf{L}\mathbf{\Lambda}_\tau^{-1}\mathbf{B}_2^{-1}
=\bra{\mathbf{L}\mathbf{\Lambda}_\tau^{-1}\mathbf{B}_2^{-1}}^T
\mathbf{L}\mathbf{\Lambda}_\tau^{-1}\mathbf{B}_2^{-1},
\end{align*}
and then
\begin{align*}
\myvec{w}_e^{T}\mathbf{\Theta}\myvec{w}_e
=\timenorm{\mathbf{L}\mathbf{\Lambda}_\tau^{-1}\mathbf{B}_2^{-1}\myvec{w}_e}^2.
\end{align*}
We use the definition \eqref{def:matrix tilde B} and the equality \eqref{decom:mateix B} to derive that
\begin{align*}
\myvec{z}_e^{T}\mathbf{Q}^{T}\mathbf{B}^{-1}\mathbf{Q}\myvec{w}_e
=&\,\bra{\brab{\mathbf{L}^{-1}}^T\mathbf{\Lambda}_\tau\mathbf{Q}\myvec{w}_e}^T
\brab{\mathbf{L}^{-1}}^T\mathbf{\Lambda}_\tau\mathbf{Q}\myvec{w}_e
=\timenorm{\brab{\mathbf{L}^{-1}}^T\mathbf{\Lambda}_\tau\mathbf{Q}\myvec{w}_e}^2\nonumber\\
=&\,\timenorm{\brab{\mathbf{L}^{-1}}^T\mathbf{\Lambda}_\tau\mathbf{Q}\mathbf{B}_2\mathbf{\Lambda}_\tau
\mathbf{L}^{-1}\mathbf{L}\mathbf{\Lambda}_\tau^{-1}\mathbf{B}_2^{-1}\myvec{w}_e}^2\nonumber\\
\le&\,\timenorm{\brab{\mathbf{L}^{-1}}^T\mathbf{\Lambda}_\tau\mathbf{Q}\mathbf{B}_2
\mathbf{\Lambda}_\tau\mathbf{L}^{-1}}^2
\timenorm{\mathbf{L}\mathbf{\Lambda}_\tau^{-1}\mathbf{B}_2^{-1}\myvec{w}_e}^2\\
=&\,\timenorm{\brab{\mathbf{L}^{-1}}^T\mathbf{Q}\widetilde{\mathbf{B}}_2\mathbf{L}^{-1}}^2\cdot
\myvec{w}_e^T\mathbf{\Theta}\myvec{w}_e\\
\le&\,\timenorm{\mathbf{Q}}\timenorm{\widetilde{\mathbf{B}}_2}\timenorm{\mathbf{L}^{-1}}^4\cdot
\myvec{w}_e^T\mathbf{\Theta}\myvec{w}_e
\le2\mathcal{M}_r\sum_{k=1}^{n}\sum_{j=1}^{k}
\theta_{k-j}^{(k)}\brat{\myvec{w}^k}^T\myvec{w}^j,
\end{align*}
where the estimate \eqref{lemproof: Quadratic inequality1} and
the definition \eqref{def: quantity Mrn} of $\mathcal{M}_r$ have been used in the last inequality.
Inserting the above inequality into \eqref{lemProof: DOC Quadratic inequality-nonlinear},
we obtain the claimed first inequality.
The second result then follows immediately
by setting $\myvec{w}^j=\myvec{z}^j$ and $\epsilon:=\sqrt{\mathcal{M}_r}$.
\end{proof}

\section{$L^2$ stability and convergence analysis}
\setcounter{equation}{0}
In this section, we shall show that $L^2$ stability and convergence analysis of the variable-step BDF2 scheme for the MBE model.
Always, they need a discrete Gr\"{o}nwall inequality \cite[Lemma 3.1]{LiaoZhang:2019linear}.
\begin{lemma}\label{lem: discrete Gronwall}
Let $\lambda\ge0$, the time sequences $\{\xi_k\}_{k=0}^N$ and $\{V_k\}_{k=1}^{N}$ be nonnegative. If
$$V_n\le \lambda\sum_{j=1}^{n-1}\tau_jV_j+\sum_{j=0}^{n}\xi_j\quad\text{for $1\leq n\le N$},$$
then it holds that
\begin{align*}
V_n\leq \exp(\lambda t_{n-1})\sum_{j=0}^{n}\xi_j
\quad\text{for\;\; $1\le n\le N$.}
\end{align*}
\end{lemma}

\subsection{$L^2$ norm stability}
We first show the $L^2$ the stability. In what follows, for notation simplicity, we shall set
$$\sum_{k, j}:=\sum_{k=1}^{n}\sum_{j=1}^{k}.$$
\begin{theorem}\label{thm: L2 norm stability}
If \textbf{S1} holds with the time-step condition $\tau_n\le \varepsilon/(16\mathcal{M}_r^2)$,
the variable-step BDF2 scheme \eqref{scheme: BDF2 implicit MBE-no-slope} is stable in the $L^2$ norm
with respect to small initial disturbance, namely,
\begin{align*}
\mynormb{\bar{\phi}^n-\phi^n}\le2\exp\brab{16\mathcal{M}_r^2t_{n-1}/\varepsilon}
\mynormb{\bar{\phi}^0-\phi^0}\quad\text{for $1\le n\le N,$}
\end{align*}
where $\bar{\phi}_h^n$ solves the equation \eqref{scheme: BDF2 implicit MBE-no-slope} with the initial data $\bar{\phi}_h^0$.
\end{theorem}
\begin{proof}Let $z_h^k$ be the solution perturbation $z_h^k:=\bar{\phi}_h^k-\phi_h^k$
for $\myvec{x}_{h}\in\bar{\Omega}_h$ and $0\le k\le N$.
Then it is easy to obtain the perturbed equation
\begin{align}\label{thmproof: L2 norm stability1}
D_{2}z_h^j+\varepsilon\Delta_{h}^2z_h^j
+\nabla_{h}\cdot\brab{\myvec{f}(\nabla_{h}\bar{\phi}_h^j)-\myvec{f}(\nabla_{h}\phi_h^j)}=0
\quad \text{for $\myvec{x}_{h}\in\Omega_h$ and $1\le j\le N$}.
\end{align}
Multiplying both sides of \eqref{thmproof: L2 norm stability1} by the DOC kernels $\theta_{k-j}^{(k)}$,
and summing up from 1 to $k$, we have
 \begin{align*}
\diff z_h^k+\varepsilon\sum_{j=1}^{k}\theta_{k-j}^{(k)}\Delta_{h}^2z_h^j
+\sum_{j=1}^{k}\theta_{k-j}^{(k)}\nabla_{h}\cdot\kbrab{\myvec{f}(\nabla_{h}\bar{\phi}_h^j)-\myvec{f}(\nabla_{h}\phi_h^j)}=0,
\end{align*}
where the equality \eqref{eq: orthogonal equality for BDF2} has been used in the derivation.
Now by taking the inner product of the above equality with $2z^k$, and summing up the derived equality from $k = 1$ to $n$, one obtain
\begin{align}\label{thmproof: L2 norm stability2}
\mynormb{z^{n}}^2-\mynormb{z^{0}}^2+2\varepsilon\sum_{k, j}\theta_{k-j}^{(k)}\myinnerb{\Delta_{h}z^j, \Delta_{h}z^k}
\le2\sum_{k, j}\theta_{k-j}^{(k)}\myinnerb{\myvec{f}(\nabla_{h}\bar{\phi}^j)-\myvec{f}(\nabla_{h} \phi^j), \nabla_{h}z^k}.
\end{align}
Now, by taking  $\myvec{v}^k:=\nabla_{h} \phi^j$ and $\myvec{z}^k:=\nabla_{h}z^k$
in the second inequality of Lemma \ref{lem: DOC Quadratic inequality-nonlinear}, one has
\begin{align}\label{thmproof: L2 nonlinearTreatment}
\mathfrak{F}\brat{\phi^n,z^n}:=&\,2\sum_{k, j}\theta_{k-j}^{(k)}
\myinnerb{\myvec{f}(\nabla_{h}\phi^j+\nabla_{h}z^j)-\myvec{f}(\nabla_{h} \phi^j),\nabla_{h}z^k}\nonumber\\
\le&\,4\sqrt{\mathcal{M}_r}\sum_{k, j}\theta_{k-j}^{(k)} \myinnerb{\nabla_{h}z^j, \nabla_{h}z^k}
=4\sqrt{\mathcal{M}_r}\sum_{k, j}\theta_{k-j}^{(k)} \myinnerb{-\Delta_{h}z^k, z^j}.
\end{align}
Note that, Lemma \ref{lem: DOC Quadratic inequality-nonlinear} holds for the simplest case $\myvec{f}(\myvec{v}):=\myvec{v}$.
Thus one can take $z^k:=-\Delta_{h}z^k$, $w^j:=z^j$ and $\epsilon=\varepsilon/(2\sqrt{\mathcal{M}_r})$ to obtain
\begin{align}\label{thmproof: L2 nonlinearTreatment2}
\mathfrak{F}\brat{\phi^n,z^n}
\le&\,2\varepsilon\sum_{k, j}\theta_{k-j}^{(k)} \myinnerb{\Delta_{h}z^j, \Delta_{h}z^k}
+8\varepsilon^{-1}\mathcal{M}_r^2\sum_{k, j}\theta_{k-j}^{(k)} \myinnerb{z^j, z^k},
\end{align}
It follows from \eqref{thmproof: L2 norm stability2} and \eqref{thmproof: L2 nonlinearTreatment2} that
\begin{align*}
\mynormb{z^{n}}^2\le\mynormb{z^0}^2
+8\varepsilon^{-1}\mathcal{M}_r^2\sum_{k, j}\theta_{k-j}^{(k)} \myinnerb{z^j, z^k}
\le \mynormb{z^0}^2+8\varepsilon^{-1}\mathcal{M}_r^2\sum_{k, j}\theta_{k-j}^{(k)} \mynormb{z^j}\mynormb{z^k}
\end{align*}
for $1\le n\le N$. Now by choosing some integer $n_1 (0 \le n_1\le n)$ such that $\mynormb{z^{n_1}}=\max_{0\le k\le n}\mynormb{z^k}$, and
setting $n= n_1 $ in the above inequality, we obtain by using Lemma \ref{lem: DOC property} (III):
\begin{align}\label{thmproof: L2 norm stability final}
\mynormb{z^{n}}\le \mynormb{z^{n_1}}\le
\mynormb{z^0}+8\varepsilon^{-1}\mathcal{M}_r^2\sum_{k=1}^{n_1}\tau_k\mynormb{z^k}
\le\mynormb{z^0}+8\varepsilon^{-1}\mathcal{M}_r^2\sum_{k=1}^{n}\tau_k\mynormb{z^k}
\end{align}
for $1\le n\le N$. By noticing the time-step condition  $\tau_n\le \varepsilon/(16\mathcal{M}_r^2)$,
one gets from \eqref{thmproof: L2 norm stability final} that
\begin{align*}
\mynormb{z^{n}}\le2\mynormb{z^0}+16\varepsilon^{-1}\mathcal{M}_r^2\sum_{k=1}^{n-1}\tau_k\mynormb{z^k}\quad\text{for $1\le n\le N$.}
\end{align*}
Then the desired result follows by using the discrete Gr\"onwall inequality in Lemma \ref{lem: discrete Gronwall}.
\end{proof}

As noticed, Theorem \ref{thm: L2 norm stability} does not involve any undesirable unbounded factors,
such as $C_r$ or $\Gamma_n$ in existing works \cite{Becker:1998,ChenWangYanZhang:2019,Emmrich:2005}.
For the time $t_n\le T$, the stability factor $\exp\brab{16\mathcal{M}_r^2t_{n-1}/\varepsilon}$ remains bounded
as the time steps $\tau_n$ vanish or the step-ratios $r_n$ approach the zero-stability limit $r_s=3.561$.
Thus Theorem \ref{thm: L2 norm stability} also shows that the variable-step BDF2 time-stepping scheme is robustly stable
with respect to the variation of time-step sizes.
Now by taking $\bar{\phi}_h^0=0$ in Theorem \ref{thm: L2 norm stability}
and using together Theorem \ref{thm: uniquely solvable},
we have the following corollary which simulates the $L^2$ norm estimate \eqref{cont: L2 norm solution estimate}.

\begin{corollary}\label{corol: L2 norm solution estimate}
If \textbf{S1} holds with the time-step condition $\tau_n\le \delta\varepsilon/(8\mathcal{M}_r^2)$ for any $0<\delta<1$,
the solution of variable-step BDF2 time-stepping scheme \eqref{scheme: BDF2 implicit MBE-no-slope} fulfills
\begin{align*}
\mynormb{\phi^n}\le\frac{1}{1-\delta}\exp\braB{\frac{8\mathcal{M}_r^2t_{n-1}}{(1-\delta)\varepsilon}}\mynormb{\phi^0}
\quad\text{for $1\le n\le N$ and $\tau_n\le \varepsilon$.}
\end{align*}
\end{corollary}

\subsection{$L^2$ norm error estimates}
We are now at the stage to give the error estimates of the variable-step BDF2 scheme. To do this, let $\xi^j:=D_2\Phi(t_j)-\partial_t\Phi(t_j)$ be the local consistency error of the BDF2 scheme at the time $t=t_j$. We will consider a convolutional consistency error $\Xi^k$ defined by
\begin{align}\label{def:BDF2-global consistency}
\Xi^k:=\sum_{j=1}^k\theta_{k-j}^{(k)}\xi^j=
\sum_{j=1}^k\theta_{k-j}^{(k)}\bra{D_2\Phi(t_j)-\partial_t\Phi(t_j)}\quad\text{for $k\ge1$.}
\end{align}
\begin{lemma}\cite[Lemma 3.4]{LiaoJiZhang:2020}\label{lem: BDF2-Consistency-Error}
If \textbf{S1} holds, the consistency error $\Xi^k$ in \eqref{def:BDF2-global consistency} satisfies
\begin{align*}
\absb{\Xi^k}
\le&\, \theta_{k-1}^{(k)}\int_{0}^{t_1} \absb{\Phi''(t)} \zd{t}
+3\sum_{j=1}^{k}\theta_{k-j}^{(k)}\tau_{j}\int_{t_{j-1}}^{t_j} \absb{\Phi'''(t)} \zd{t}\quad\text{for $k\ge1$}
\end{align*}
such that
\begin{align*}
\sum_{k=1}^n\absb{\Xi^k}
\le&\, \tau_1\int_{0}^{t_1} \absb{\Phi''(t)} \zd{t}\,\sum_{k=1}^n\prod_{i=2}^k\frac{r_i^2}{1+2r_i}
+3t_n\max_{1\le j\le n}\braB{\tau_{j}\int_{t_{j-1}}^{t_j} \absb{\Phi'''(t)} \zd{t}}\quad\text{for $n\ge1$.}
\end{align*}
\end{lemma}

Hereafter, we shall use a generic constant $C_{\phi}>0$ in the error estimates which is
not necessarily the same at different occurrences, but always independent of the time steps $\tau_n$, the step-ratios $r_n$ and the spatial length $h$.

\begin{theorem}\label{thm: L2 norm convergence}
Assume that the MBE problem \eqref{cont: MBE no-slope problem}
has a smooth solution $\Phi\in C_{\myvec{x},t}^{(6,3)}(\Omega\times(0,T])$.
If \textbf{S1} holds with the time-steps $\tau_n\le \varepsilon/(16\mathcal{M}_r^2)$,
the BDF2 scheme \eqref{scheme: BDF2 implicit MBE-no-slope}
admits the following error estimate:
\begin{align*}
\mynormb{\Phi^n-\phi^n}\le&\,C_{\phi}\exp(16\mathcal{M}_r^2t_{n-1}/\varepsilon)
\kbra{\tau_1^2\sum_{k=1}^n\prod_{i=2}^k\frac{r_i^2}{1+2r_i}+t_n(\tau^2+h^2)}\quad\text{for $1\le n\le N$.}
\end{align*}
\end{theorem}
\begin{proof}
Let $\Phi_{h}^n:=\Phi(\myvec{x}_h,t_n)$ and $e_{h}^n$  be the error function
$e_{h}^n:=\Phi_{h}^n-\phi_{h}^n$ with $e_{h}^n:=0$
for $\myvec{x}_h\in\bar{\Omega}_h$.
We then have the following error equation
\begin{align}\label{thmproof: error equation-1}
D_{2}e_h^j+\varepsilon\Delta_{h}^2e_h^j
+\nabla_{h}\cdot\kbrab{\myvec{f}(\nabla_{h}\Phi_h^j)-\myvec{f}(\nabla_{h}\phi_h^j)}
=\xi_h^j+\eta_h^j,
\end{align}
where $\xi_h^j$ and $\eta_h^j$ are the local consistency error in time and physical domain, respectively. If the solution is smooth,
Lemma \ref{lem: DOC property} (III) gives
\begin{align}\label{thmproof: consistency-space}
\sum_{k=1}^n\mynormb{\Pi^k}\le C_{\phi}t_nh^2\quad\text{for $1\le n\le N$,\quad where $\Pi_h^k:=\sum_{j=1}^k\theta_{k-j}^{(k)}\eta_h^j$.}
\end{align}
Multiplying both sides of \eqref{thmproof: error equation-1} by the DOC kernels $\theta_{k-j}^{(k)}$,
and summing up the superscript from $j=1$ to $k$,  we obtain by applying the equality \eqref{eq: orthogonal equality for BDF2}
\begin{align*}
\diff e_h^k+\varepsilon\sum_{j=1}^{k}\theta_{k-j}^{(k)}\Delta_{h}^2e_h^j
+\sum_{j=1}^{k}\theta_{k-j}^{(k)}\nabla_{h}\cdot
\kbrab{\myvec{f}(\nabla_{h} \Phi_h^j)-\myvec{f}(\nabla_{h}\phi_h^j)}=\Xi_h^k+\Pi_h^k,
\end{align*}
where $\Xi_h^k$ and $S_h^k$ are defined by \eqref{def:BDF2-global consistency} and \eqref{thmproof: consistency-space}, respectively.
Now by taking the inner product of the above equality with $2e^k$,
and summing up the superscript from $k=1$ to $n$, we obtain by using the discrete Green's formula
\begin{align}\label{thmproof: error equation-3}
\mynormb{e^n}^2-\mynormb{e^0}^2+2\varepsilon\sum_{k,j}\theta_{k-j}^{(k)}
\myinnerb{\Delta_{h}e^j,\Delta_{h}e^k}\le\mathfrak{F}\brat{\phi^n,e^n}+2\sum_{k=1}^{n}\myinnerb{\Xi^k+\Pi^k, e^k},
\end{align}
where $\mathfrak{F}\brat{\phi^n,e^n}$ is defined in \eqref{thmproof: L2 nonlinearTreatment}. 
The derivation of \eqref{thmproof: L2 nonlinearTreatment2} yields
\begin{align*}
\mathfrak{F}\brat{\phi^n,e^n}\le&\,2\varepsilon\sum_{k, j}\theta_{k-j}^{(k)} \myinnerb{\Delta_{h}e^j, \Delta_{h}e^k}
+8\varepsilon^{-1}\mathcal{M}_r^2\sum_{k, j}\theta_{k-j}^{(k)} \myinnerb{e^j, e^k}.
\end{align*}
With the help of Cauchy-Schwarz inequality, it follows from \eqref{thmproof: error equation-3} that
\begin{align}\label{thmproof: error equation-4}
\mynormb{e^n}^2\le\mynormb{e^0}^2
+8\varepsilon^{-1}\mathcal{M}_r^2\sum_{k, j}\theta_{k-j}^{(k)}\mynormb{e^j}\mynormb{e^k}
+2\sum_{k=1}^{n}\mynormb{\Xi^k+\Pi^k}\mynormb{e^k}\quad\text{for $1\le n\le N$.}
\end{align}
Then, by choosing some integer $n_2 (0 \le n_2\le n)$ such that
$\mynormb{e^{n_2}}=\max_{0\le k\le n}\mynormb{e^k},$ and
setting $n= n_2$ in the above inequality \eqref{thmproof: error equation-4},
we obtain by using together Lemma \ref{lem: DOC property} (III) and the time-step condition
$\tau_n\le\varepsilon/(16\mathcal{M}_r^2)$
\begin{align*}
\mynormb{e^{n}}\le2\mynormb{e^0}
+16\varepsilon^{-1}\mathcal{M}_r^2\sum_{k=1}^{n-1}\tau_{k}\mynormb{e^k}
+4\sum_{k=1}^{n}\mynormb{\Xi^k+\Pi^k}\quad\text{for $1\le n\le N$.}
\end{align*}
Then by the discrete Gr\"onwall inequality in Lemma \ref{lem: discrete Gronwall} we have
\begin{align*}
\mynormb{e^n}\le2\exp(16\mathcal{M}_r^2t_{n-1}/\varepsilon)\braB{\mynormb{e^0}
+2\sum_{k=1}^n\mynormb{\Xi^k}+2\sum_{k=1}^n\mynormb{\Pi^k}}
\quad\text{for $1\le n\le N$.}
\end{align*}
The desired result follows by using together the estimates in \eqref{thmproof: consistency-space}
and Lemma \ref{lem: BDF2-Consistency-Error}.
\end{proof}

Notice that Theorem \ref{thm: L2 norm convergence} confirms at least a first-order convergence rate of the numerical solution under the step-ratio condition \textbf{S1}, as  $\tau_1\sum_{k=1}^n\prod_{i=2}^k\frac{r_i^2}{1+2r_i}\le t_n.$ While the second-order rate of convergence can be recovered if the following assumption is fulfilled:
  \begin{enumerate}[itemindent=1em]
\item[\textbf{S2}.]
  The time-step ratios $r_k$ are contained in  \textbf{S1},
  but almost all of them less than $1+\sqrt{2}$, or $\abs{\mathfrak{R}}=N_0\ll N$,
  where $\mathfrak{R}$ is an index set  $ \mathfrak{R}:=\{k|1+\sqrt{2}\le r_k < (3+\sqrt{17})/2\}. $
\end{enumerate}
Although the condition \textbf{S1} allows one to use a series of increasing time-steps
with the amplification factors up to 3.561, while in practice, the use of large time-steps will in general result in a loss of numerical
accuracy. In this sense, the condition \textbf{S2} is much more reasonable in practice because
large amplification factors of time-step size are rarely appeared continuously
in long-time simulations. As shown in \cite[Lemma 3.3]{LiaoZhang:2019linear},
there exists a step-ratio-dependent constant $c_r$ such that
$$\sum_{k=1}^n\prod_{i=2}^k\frac{r_i^2}{1+2r_i}\le c_r.$$
This results in the following corollary.
\begin{corollary}\label{corol: L2 norm convergence}
 Assume that the nonlinear MBE problem \eqref{cont: MBE no-slope problem} has a unique smooth solution.
 If the step-ratio assumption  \textbf{S2} holds with the time-steps $\tau_n\le\varepsilon/(16\mathcal{M}_r^2)$,
the BDF2 scheme \eqref{scheme: BDF2 implicit MBE-no-slope} is second-order convergent in the $L^2$ norm,
\begin{align*}
\mynormb{\Phi^n-\phi^n}\le&\,C_{\phi}\exp(16\mathcal{M}_r^2t_{n-1}/\varepsilon)
\braB{c_r\tau_1^2+t_n(\tau^2+h^2)}\quad\text{for $1\le n\le N$.}
\end{align*}
\end{corollary}

\section{Numerical examples}
\setcounter{equation}{0}
In this section, we shall present some numerical experiments to verify our theoretical findings. In all our computations, a fixed-point iteration scheme will be employed to solve the nonlinear BDF2 scheme at each time level with a tolerance $10^{-12}$.
\subsection{Random generated time meshes}
We first test the performance on random generated time meshes.
To this end, we set $\varepsilon = 0.1 $ and consider the following exterior-forced MBE model
$$\Phi_t=-\varepsilon\Delta^2\Phi-\nabla\cdot \myvec{f}(\nabla\Phi)+g(x,t), \quad \Omega=(0,2\pi)^2. $$
The function $g(x,t)$ is chosen such that the exact solution yields $\Phi(x,t)=\cos(t)\sin(x)\sin(y).$
The accuracy of the variable-step BDF2 scheme is tested via the random meshes. Let
$$\tau_k:= T\sigma_k/S, \quad 1\leq k\leq N,$$
where $\sigma_k\in(0,1)$ is a uniformly distributed random number and $S=\Sigma_{k=1}^{N}\sigma_k$.
The discrete error in the $L^2$-norm will be tested: $e(N):=\mynorm{\Phi(T)-\phi^N},$  and the following convergence rate will be reported:
$$ \text{Order}  \approx \log(e(N)/e(2N))/\log(\tau(N)/\tau(2N)),$$ where $\tau(N)$ denotes the maximal time-step size for total $N$ subintervals.

\begin{table}[htb!]
\begin{center}\label{MBE-BDF2-Time-Error}
\caption{Accuracy of BDF2 scheme \eqref{scheme: BDF2 implicit MBE-no-slope} on random time mesh.}\vspace*{0.3pt}
\def\temptablewidth{0.7\textwidth}
{\rule{\temptablewidth}{0.5pt}}
\begin{tabular*}{\temptablewidth}{@{\extracolsep{\fill}}cccccc}
  $N$   &$\tau$      &$e(N)$     &Order  &$\max r_k$  &$N_1$\\
  \midrule
  10    &1.49e-01    &1.23e-01   &$-$    & 2.94  &0\\
  20    &9.16e-02    &8.20e-02   &1.84   &11.98  &3\\
  40    &5.52e-02    &2.57e-02   &2.29   &34.82  &7\\
  80    &2.70e-02    &4.78e-03   &2.35   &37.72  &13\\
  160   &1.23e-02    &7.20e-04   &2.42   &71.89  &24\\
  320   &6.26e-03    &1.85e-04   &2.00   &850.80 &49\\
\end{tabular*}
{\rule{\temptablewidth}{0.5pt}}
\end{center}
\end{table}	

In this example, we use 3000 grid points in the physical domain
and solve the problem until $T = 1$. The numerical results are presented in Table \ref{MBE-BDF2-Time-Error},
in which we have also recorded the maximal time-step size $\tau$ , the
maximal step ratio and the number (denote by $N_1$ in Table 1) of time levels with the step ratio
$r_k \geq (3 + \sqrt{17})/2$. It is clear seen that the BDF2 scheme admits a
second-order rate of convergence for those nonuniform time meshes.

\subsection{Adaptive time-stepping strategy}
\begin{algorithm}
\caption{Adaptive time-stepping strategy}
\label{Adaptive-Time-Step-Strategy}
\begin{algorithmic}[1]
\Require{Given $\phi^{n}$ and time step $\tau_{n}$}
\State Compute $\phi_{1}^{n+1}$ by using BDF1 scheme with time step $\tau_{n}$.
\State Compute $\phi_{2}^{n+1}$ by using BDF2 scheme with time step $\tau_{n}$.
\State Calculate $e_{n+1}=\|\phi_{2}^{n+1}-\phi_{1}^{n+1}\|/\|\phi_{2}^{n+1}\|$.
\If {$e_{n+1}<tol$ or $\tau_{n}\le{\tau_{\min}}}$
\If {$e_{n+1}<tol$}
\State Update time-step size
 $\tau_{n+1}\leftarrow\min\{\max\{\tau_{\min},\tau_{ada}\},\tau_{\max}\}$.
\Else
\State Update time-step size $\tau_{n+1}\leftarrow\tau_{\min}$.
\EndIf
\Else
\State Recalculate with time-step size
 $\tau_{n}\leftarrow\min\{\max\{\tau_{\min},\tau_{ada}\},\tau_{\max}\}$; Goto 1.
\EndIf
\end{algorithmic}
\end{algorithm}

Next we test a practical adaptive time-stepping strategy in \cite{GomezHughes:2011Provably}.
Different  adaptive time-stepping strategies can also be found in
\cite{QiaoZhangTang:2011,ZhangMaQiao:2013}.
As verified in the previous sections, the variable-step BDF2 scheme
\eqref{scheme: BDF2 implicit MBE-no-slope} is robustly stable with respect to
the step-size variations satisfying the step-ratio condition \textbf{S1}.
In \cite{GomezHughes:2011Provably}, the adaptive time-step $\tau_{ada}$
(the next step) is updated adaptively using the current step information $\tau_{cur}$ via the following formula
\begin{align*}
\tau_{ada}\bra{e,\tau_{cur}}
=\min\big\{S_a\sqrt{tol/e}\,\tau_{cur},r_s\tau_{cur}\big\}
\end{align*}
where $e$ is the relative error of solution at the current time-level,
$tol$ is a reference tolerance, $S_a$ is some default safety parameter determined by try-and-error tests.
Notice that $r_s=3.561$ is an artificial constant that is due to the condition \textbf{S1}.
More details of the above adaptive time-stepping strategy can be found in Algorithm \ref{Adaptive-Time-Step-Strategy}.
In our computation, if not explicitly specified, we choose the safety coefficient as $S_a=0.9$,
and set the reference tolerance $tol=10^{-3}.$ The maximal time step is chosen as $\tau_{\max}=0.1$ which the minimal time step is set to be $\tau_{\min}=10^{-4}$.

In this example, we consider the MBE model (\ref{cont: MBE no-slope problem}) with the following initial condition
\begin{align}\label{example-2}
\phi_0(x,y)=0.1(\sin3x\sin2y+\sin5x\sin5y).
\end{align}
We take the parameter $\varepsilon=0.1$ and use a $128\times 128$ uniform mesh in the physical domain $\Omega=(0,2\pi)^2$.
To obtain the deviation of the height function, we define the roughness measure function $R(t)$ as follow, $R(t)=\sqrt{\frac{1}{\abs{\Omega}}\int_{\Omega}(\phi(\myvec{x},t)-\bar{\phi}(\myvec{x},t))^2\zd\myvec{x}}$, where $\bar{\phi}(t)=\frac{1}{\abs{\Omega}}\int_{\Omega}\phi(\myvec{x},t)\zd\myvec{x}$ is the average.

\begin{figure}[htb!]
\centering
\includegraphics[width=2.0in]{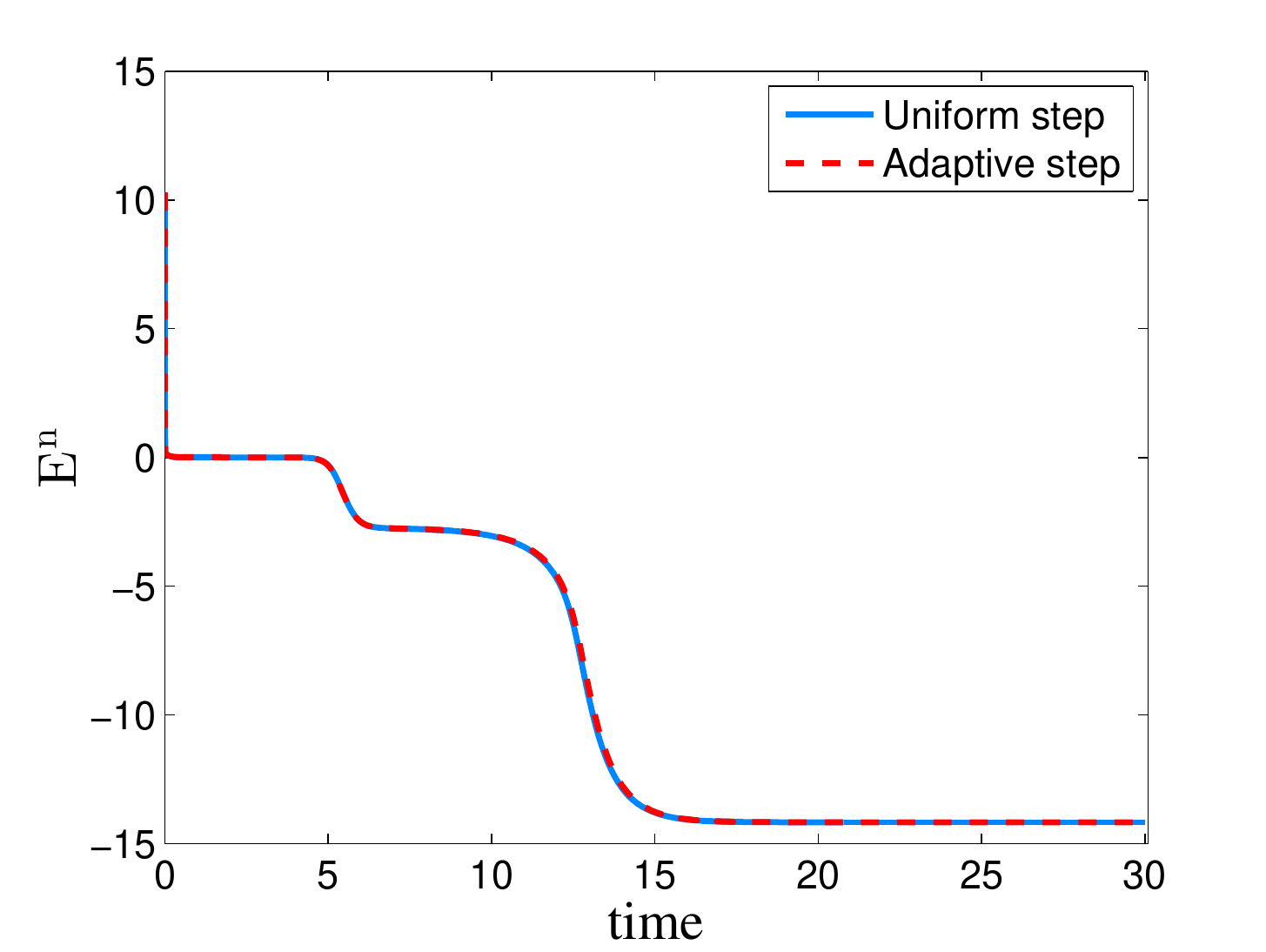}
\includegraphics[width=2.0in]{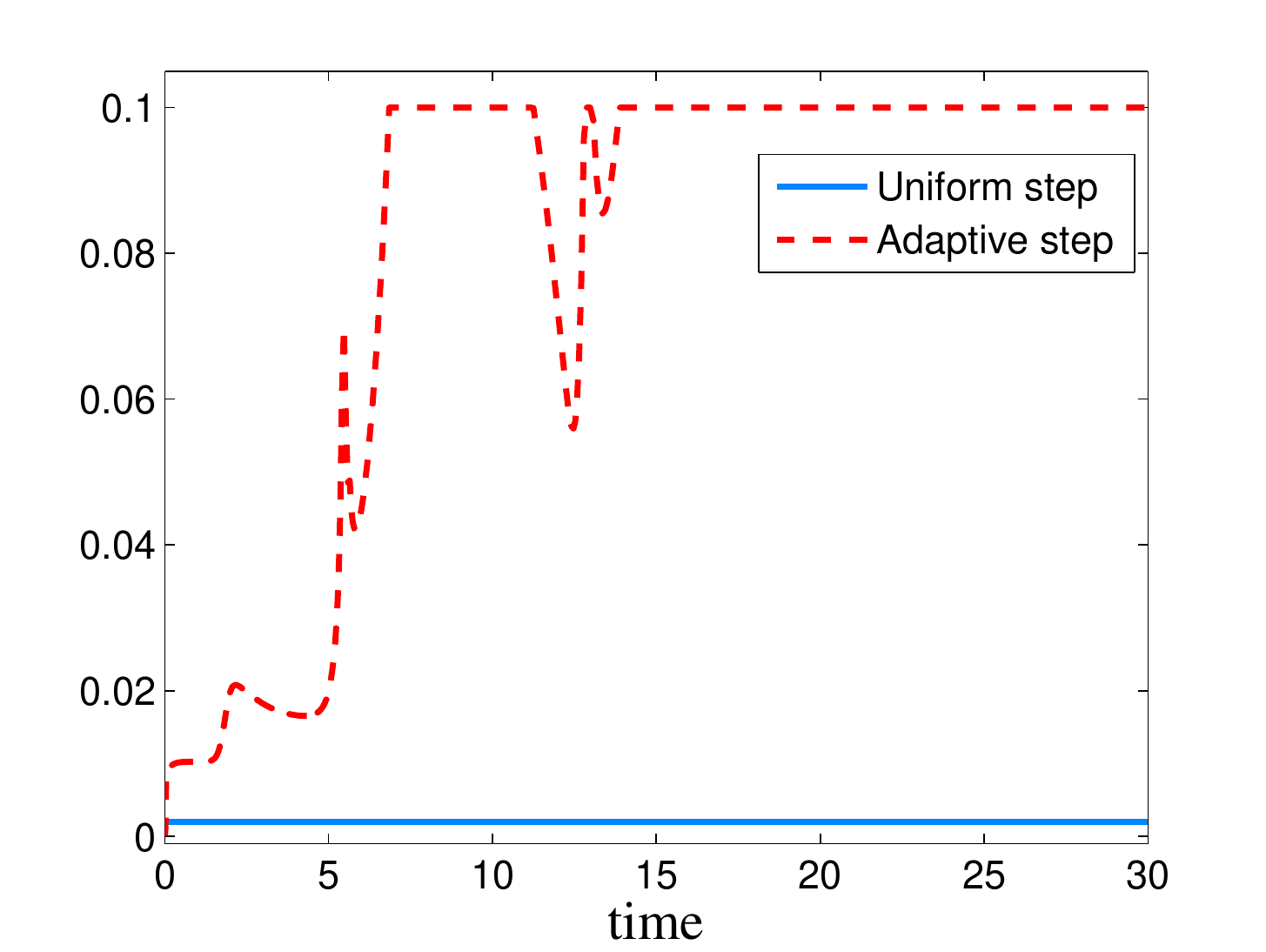}
\includegraphics[width=2.0in]{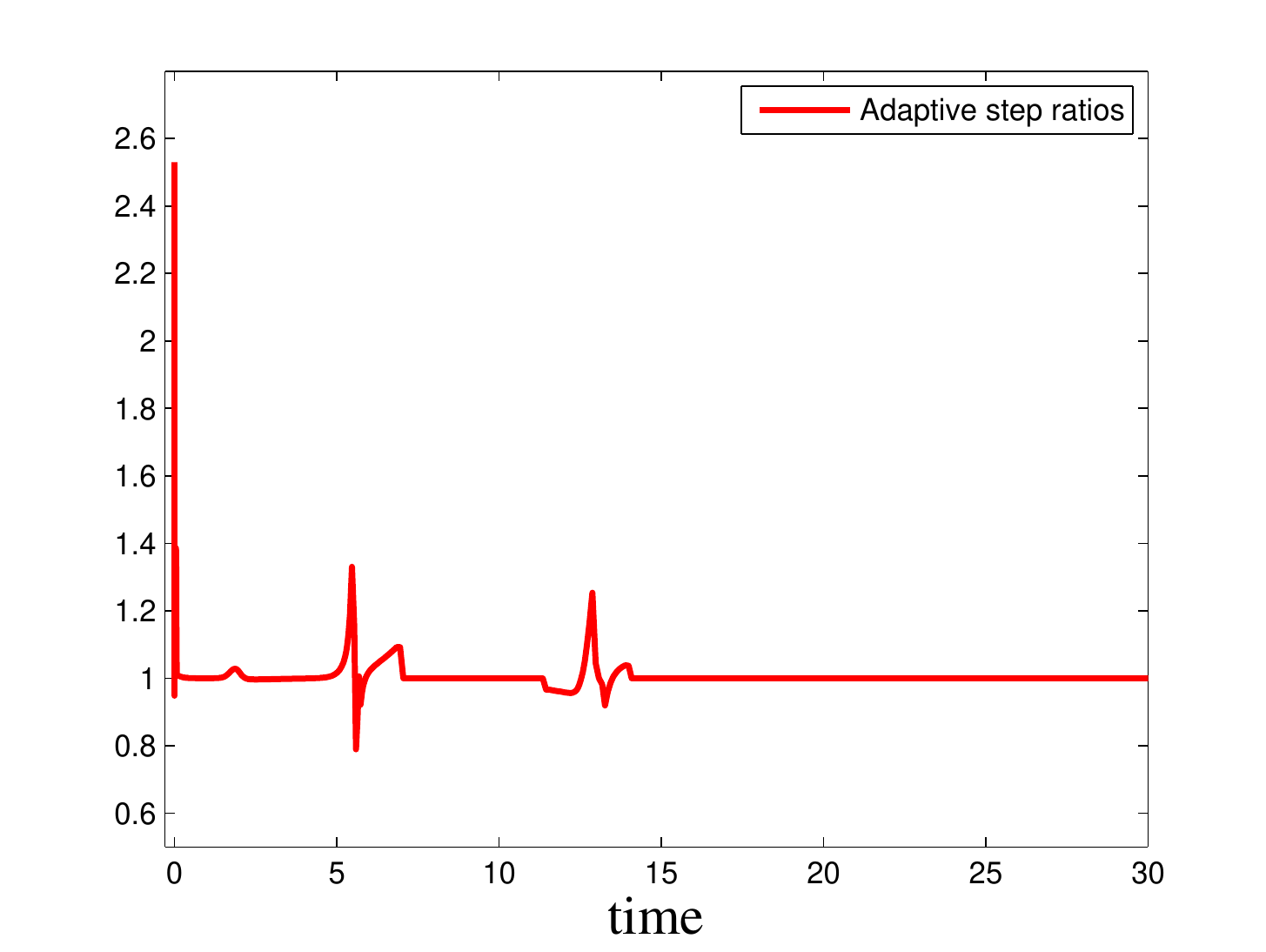}\\
\caption{Evolutions of energy (left), time steps (middle)  and time steps ratios (right) 
of  the MBE equation using different time strategies until time T=30.}
\label{Comparison-Uniform-Adaptive-Energy}
\end{figure}

We aim at simulating the benchmark problem with an initial condition of \eqref{example-2}.
We first test the efficiency and accuracy of Algorithm \ref{Adaptive-Time-Step-Strategy}.
To make a comparison, we shall also show the numerical
results with the uniform time meshes.
The solution is first simulated until $T=30$ with a constant time step $\tau=10^{-3}$.
We then use the adaptive time-stepping strategy described in Algorithm 1 to
repeat the simulation.
The numerical results are summarized in Figure \ref{Comparison-Uniform-Adaptive-Energy}.
We note that it takes 30000 uniform time steps with $\tau = 10^{-3}$,
while the total number of adaptive time steps is only 529 to get the similar results,
meaning that the time-stepping adaptive strategy is computationally efficient.
In addition, the right subplot in Figure \ref{Comparison-Uniform-Adaptive-Energy}
shows that the adaptive step-ratios satisfy the condition $\textbf{S1}$.

\begin{figure}[htb!]
\centering
\includegraphics[width=2.0in]{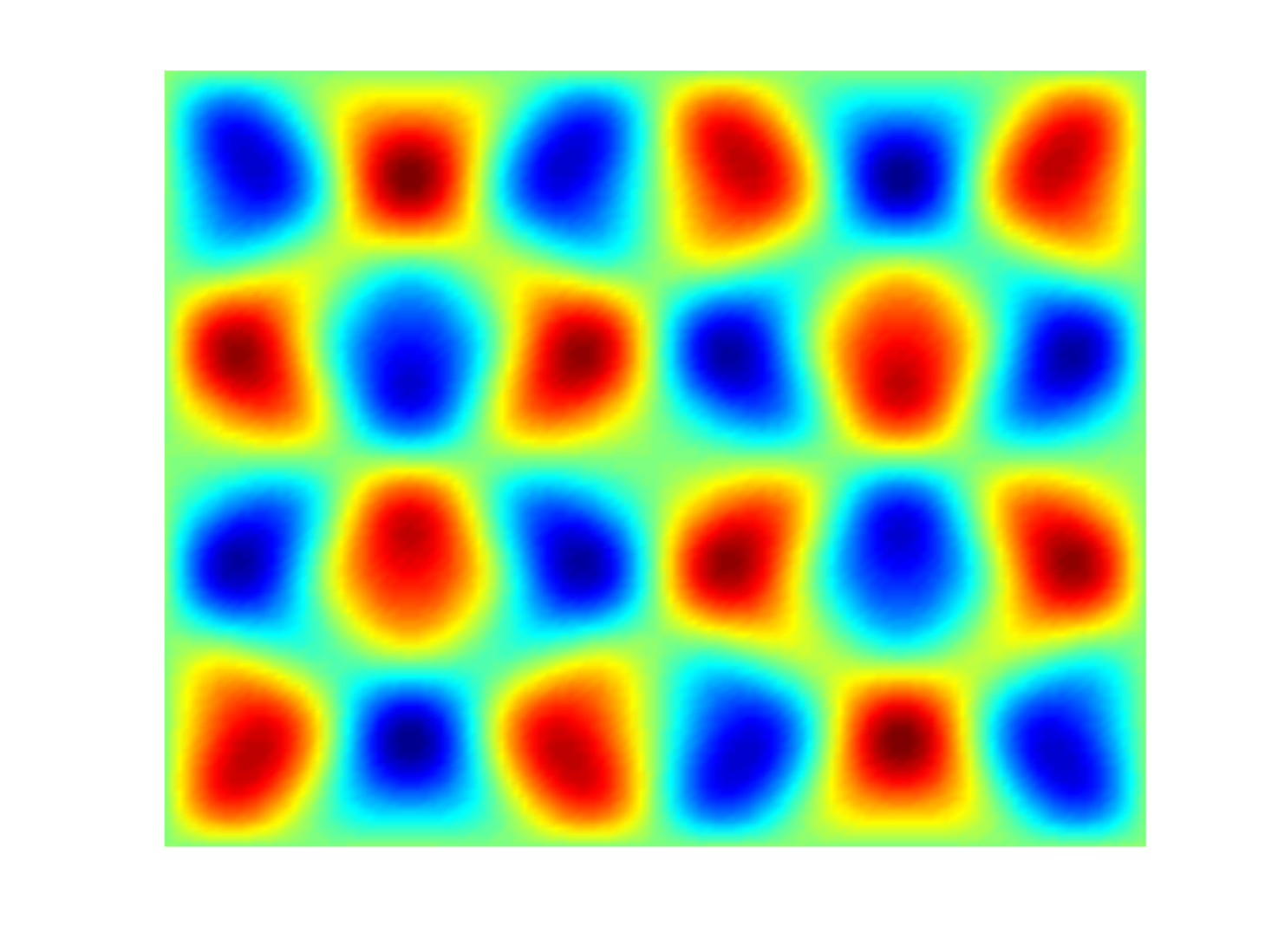}
\includegraphics[width=2.0in]{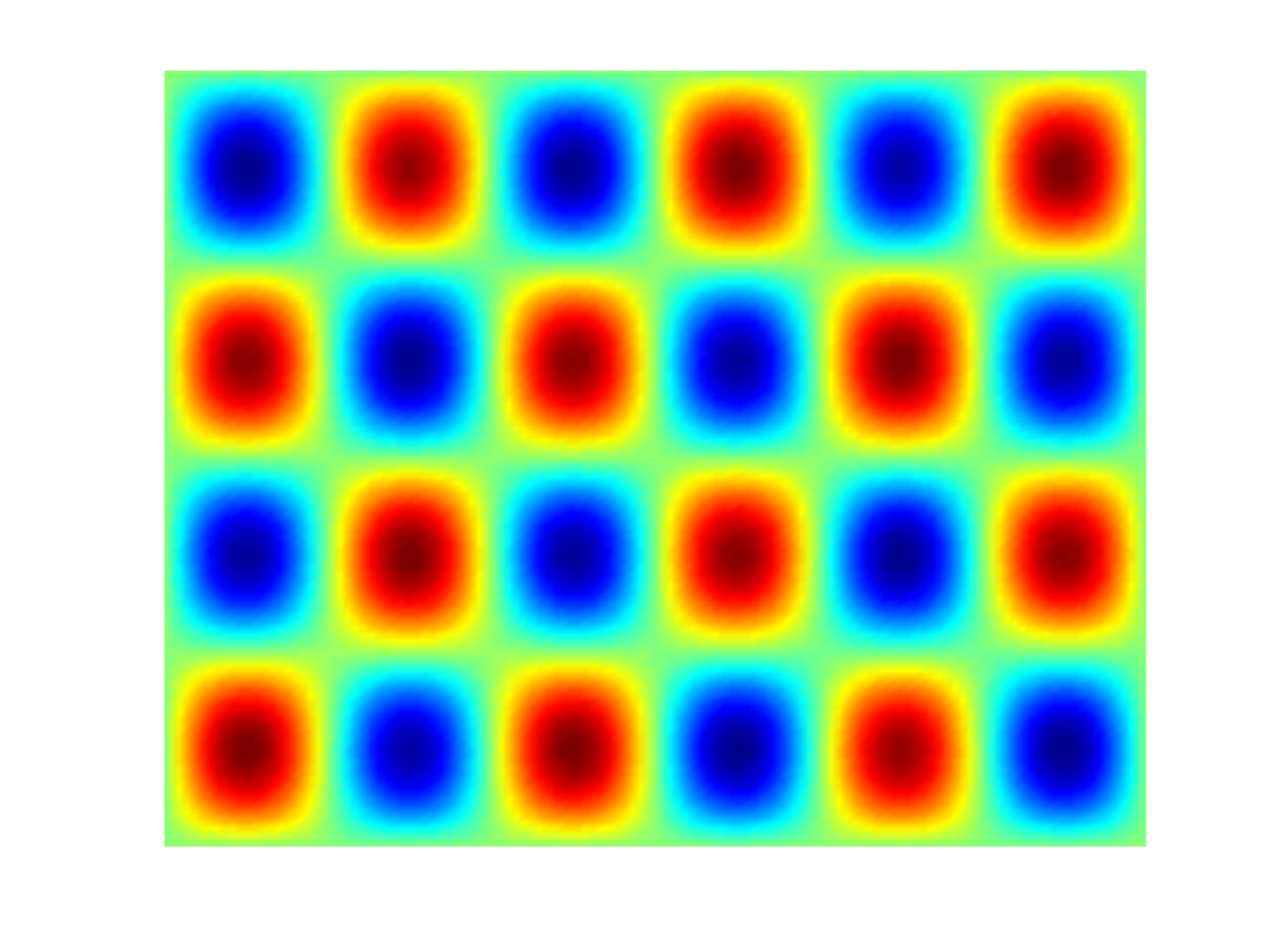}
\includegraphics[width=2.0in]{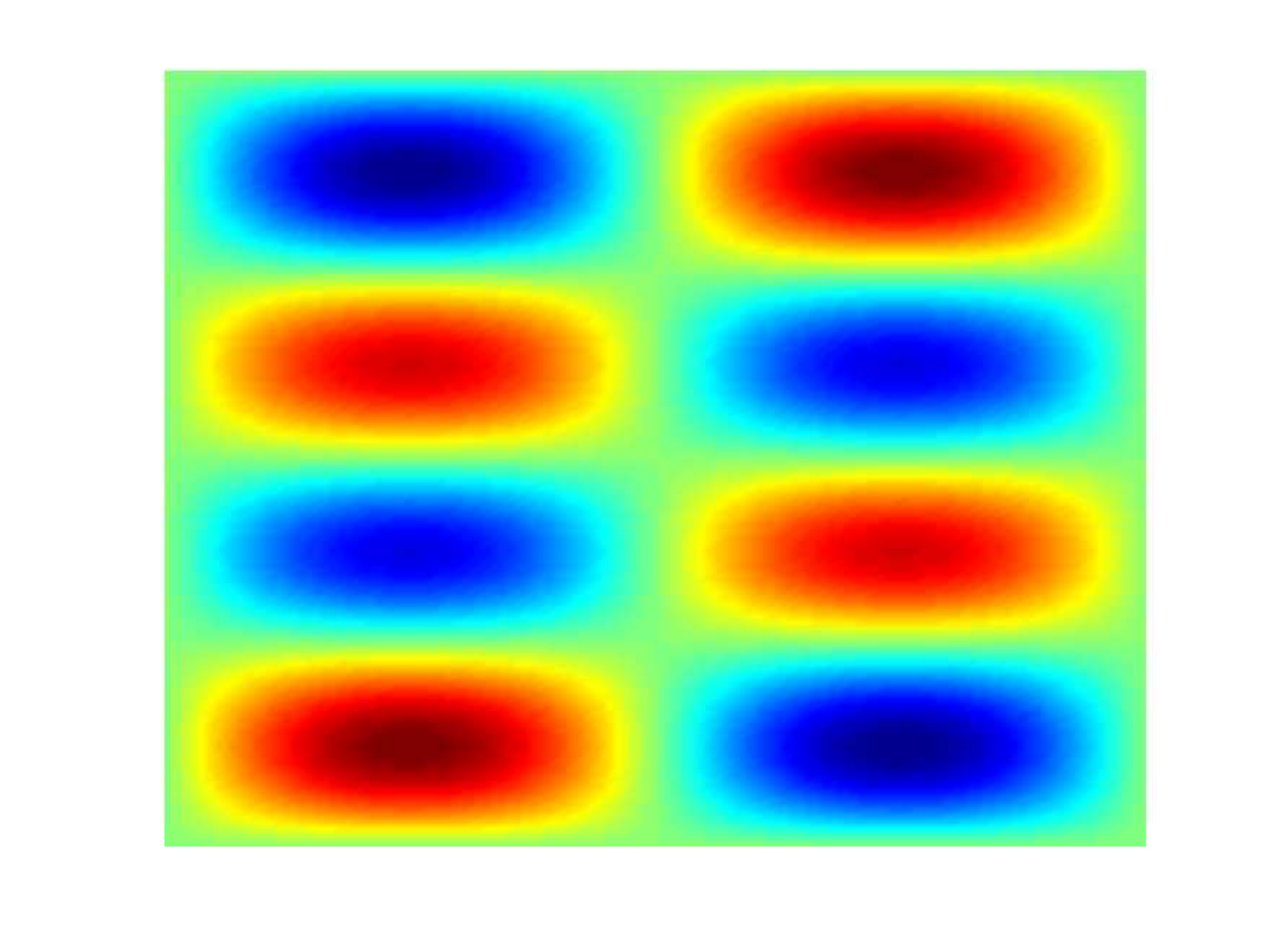}\\
\includegraphics[width=2.0in]{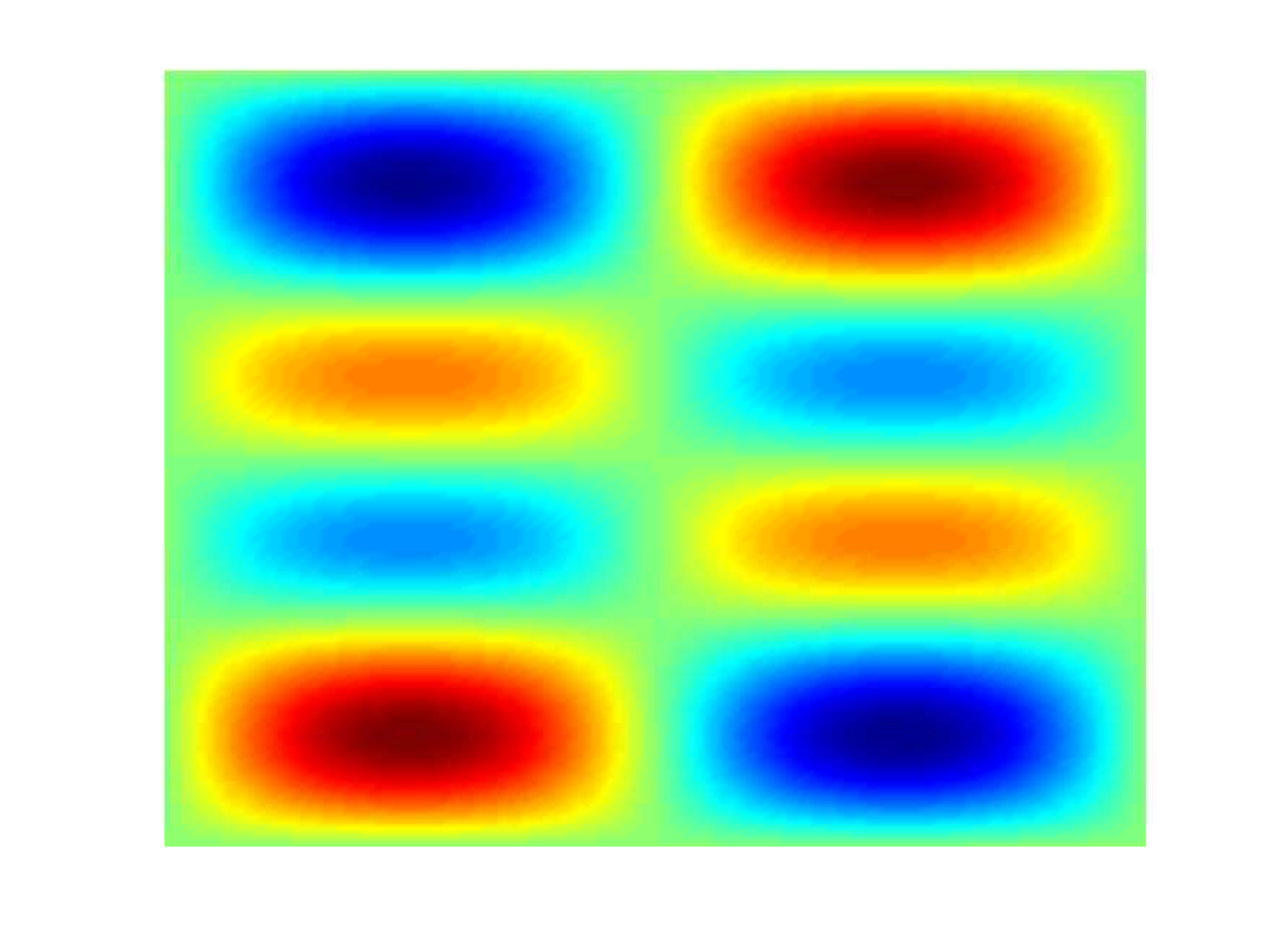}
\includegraphics[width=2.0in]{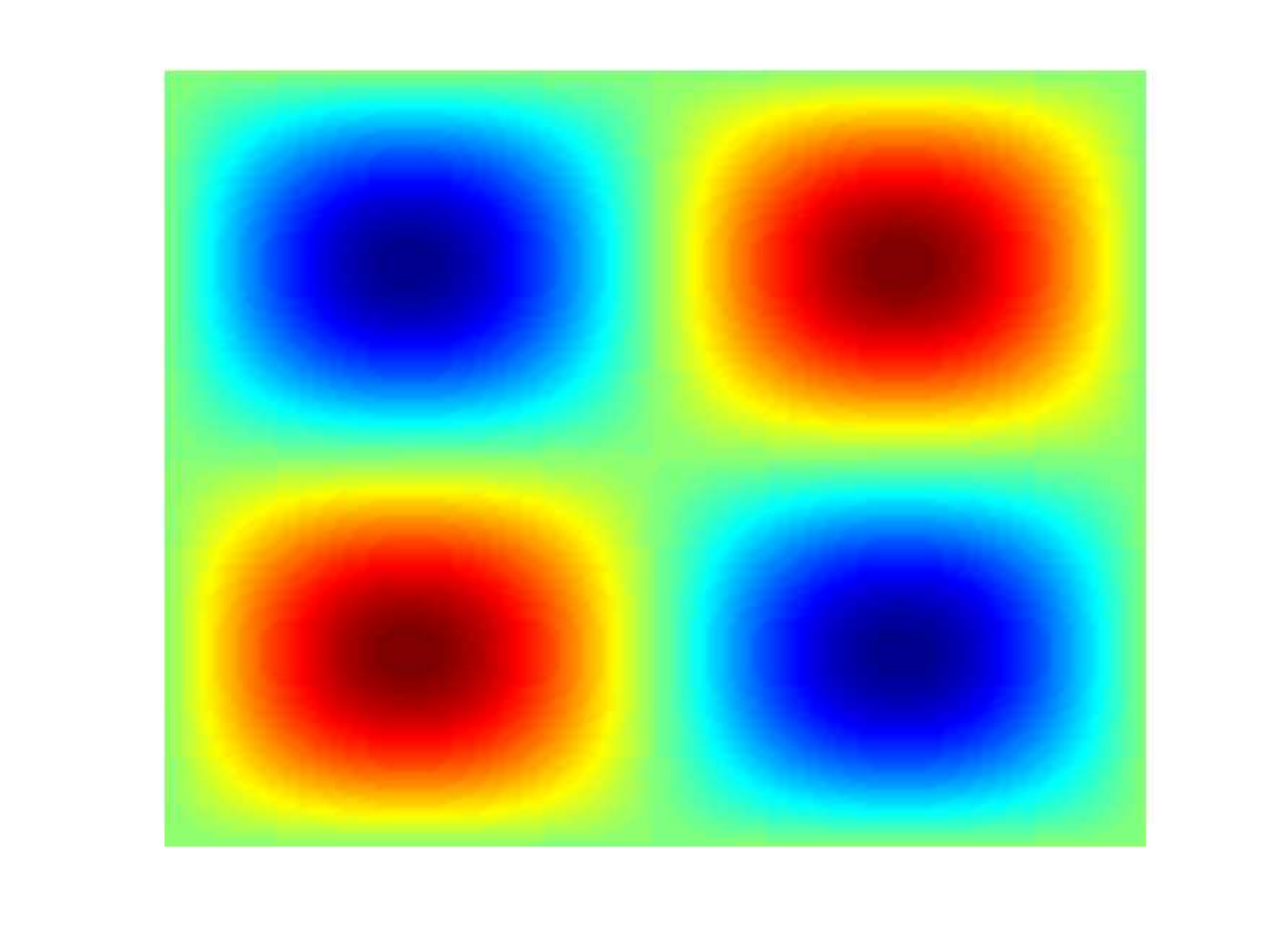}
\includegraphics[width=2.0in]{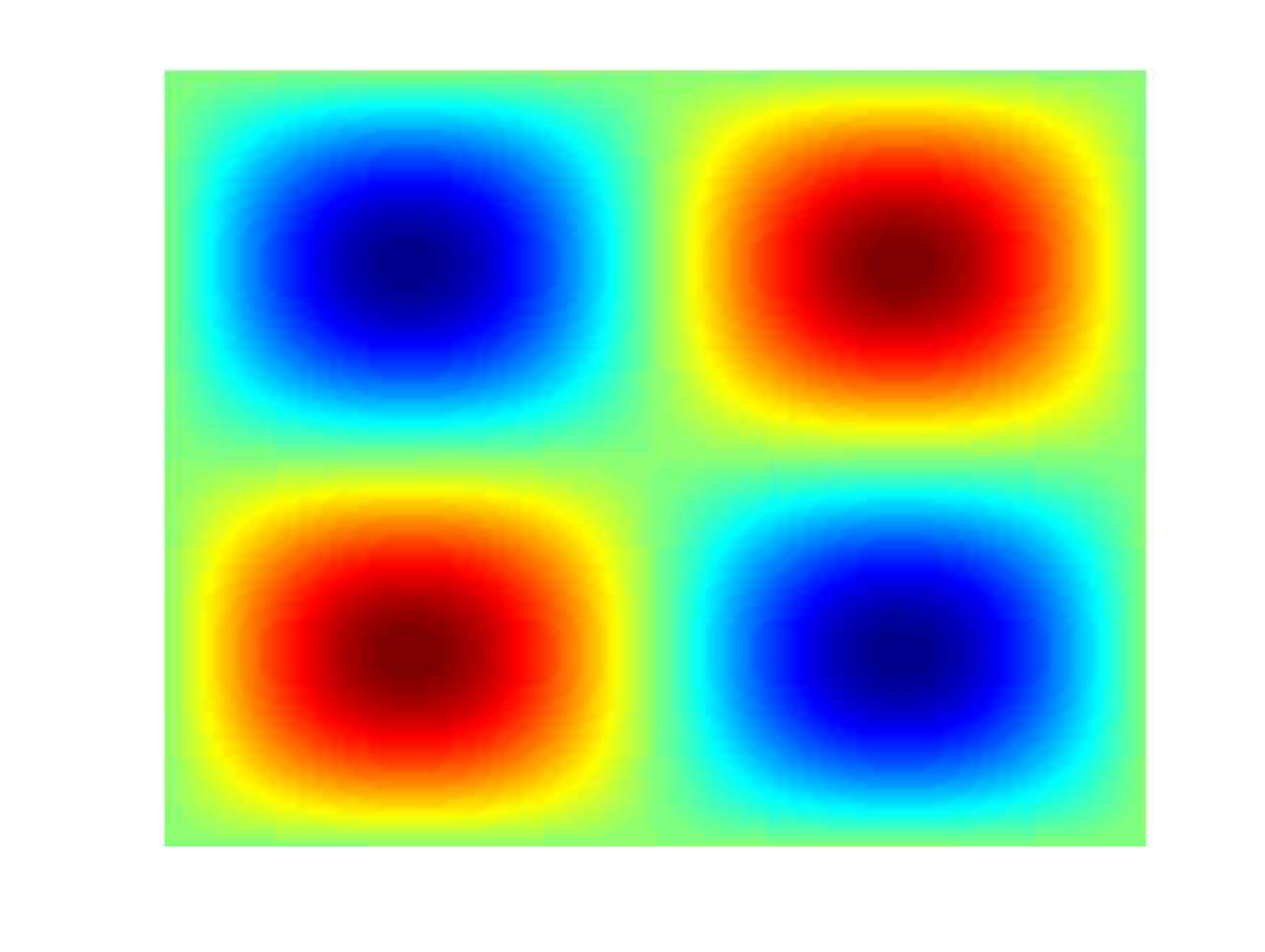}\\
\caption{The isolines of numerical solutions of the height function $\phi$ for
   the MBE equation using adaptive time strategy at $t=0, 1, 5, 10, 20,30$, respectively.}
\label{example2-Snapshots-Dynamics}
\end{figure}

\begin{figure}[htb!]
\centering
\includegraphics[width=2.0in]{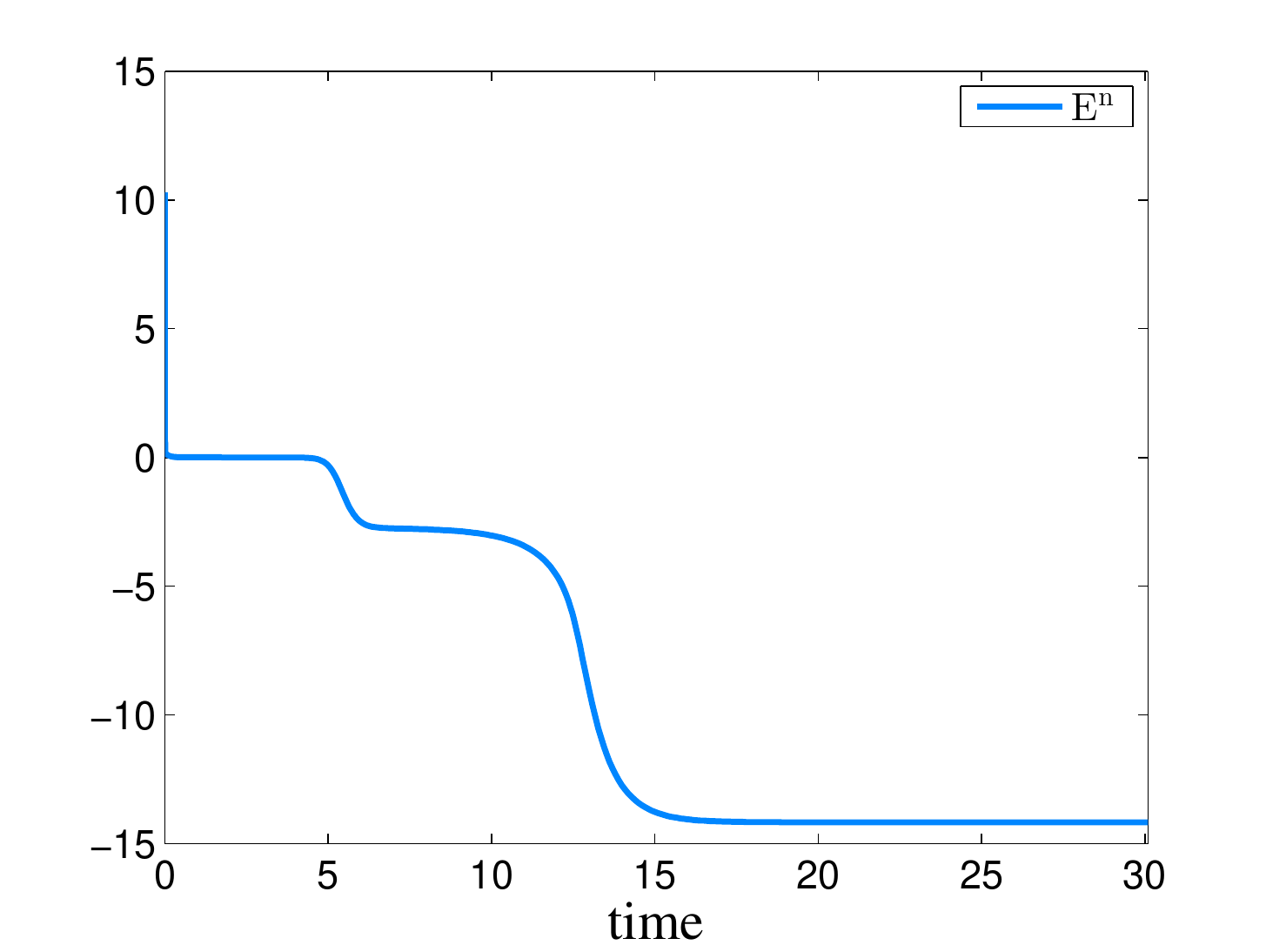}
\includegraphics[width=2.0in]{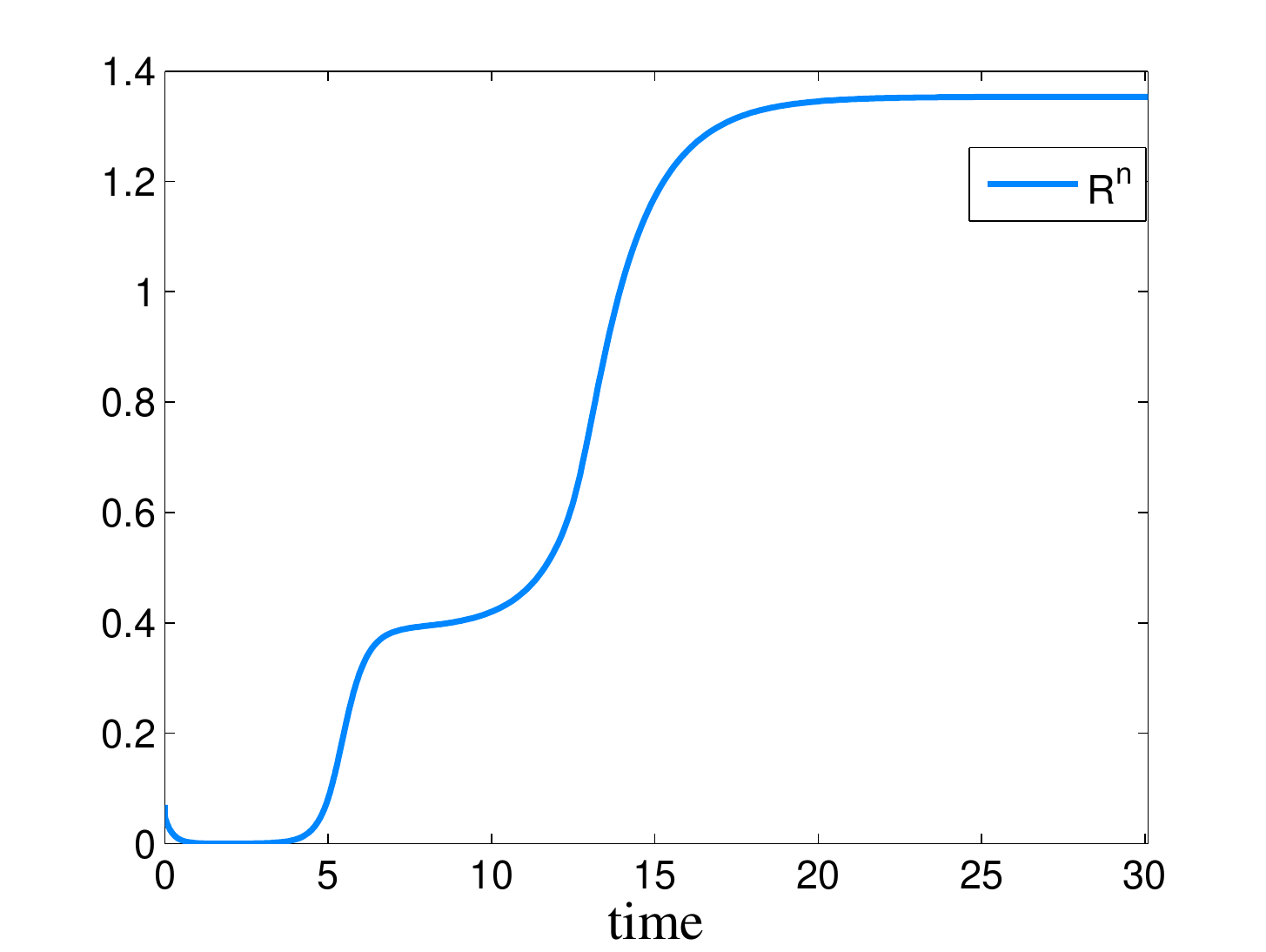}
\includegraphics[width=2.0in]{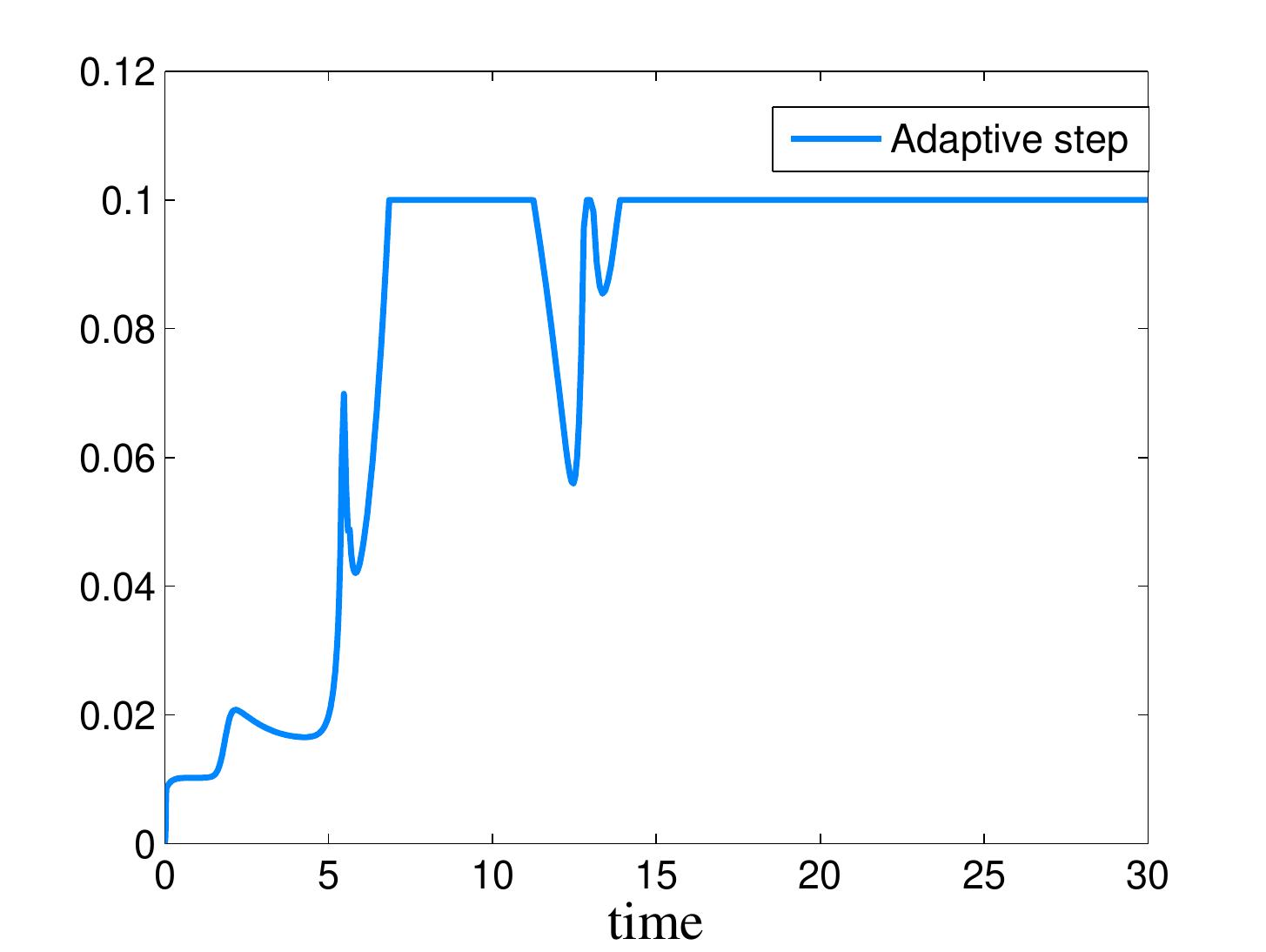}\\
\caption{Evolutions of energy (left), roughness (middle) and adaptive
  time steps (right) for the MBE equation using adaptive time strategy, respectively}
\label{example2-Adaptive-Energy-Roughness}
\end{figure}

\begin{figure}[htb!]
\centering
\includegraphics[width=3.0in]{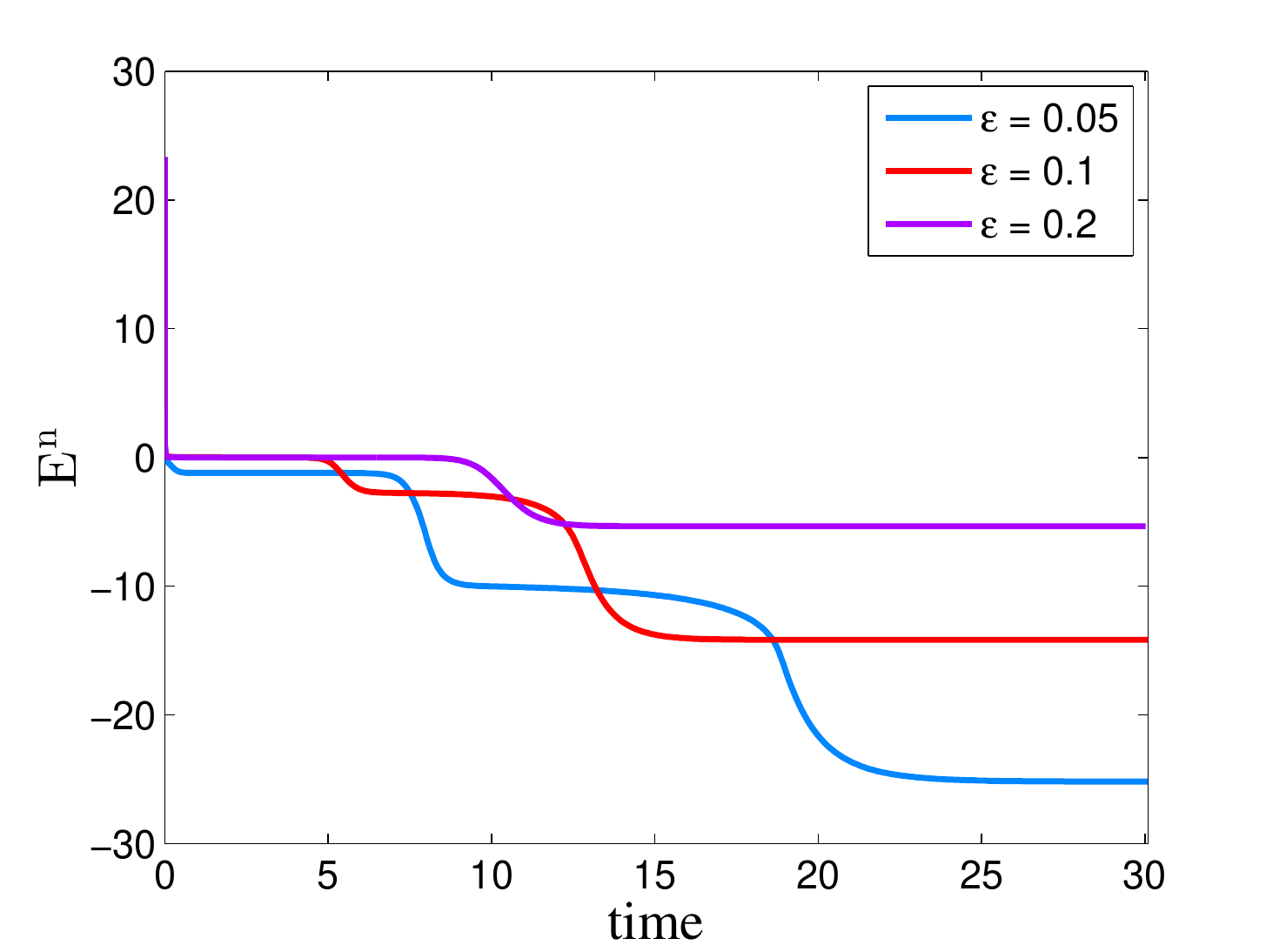}
\includegraphics[width=3.0in]{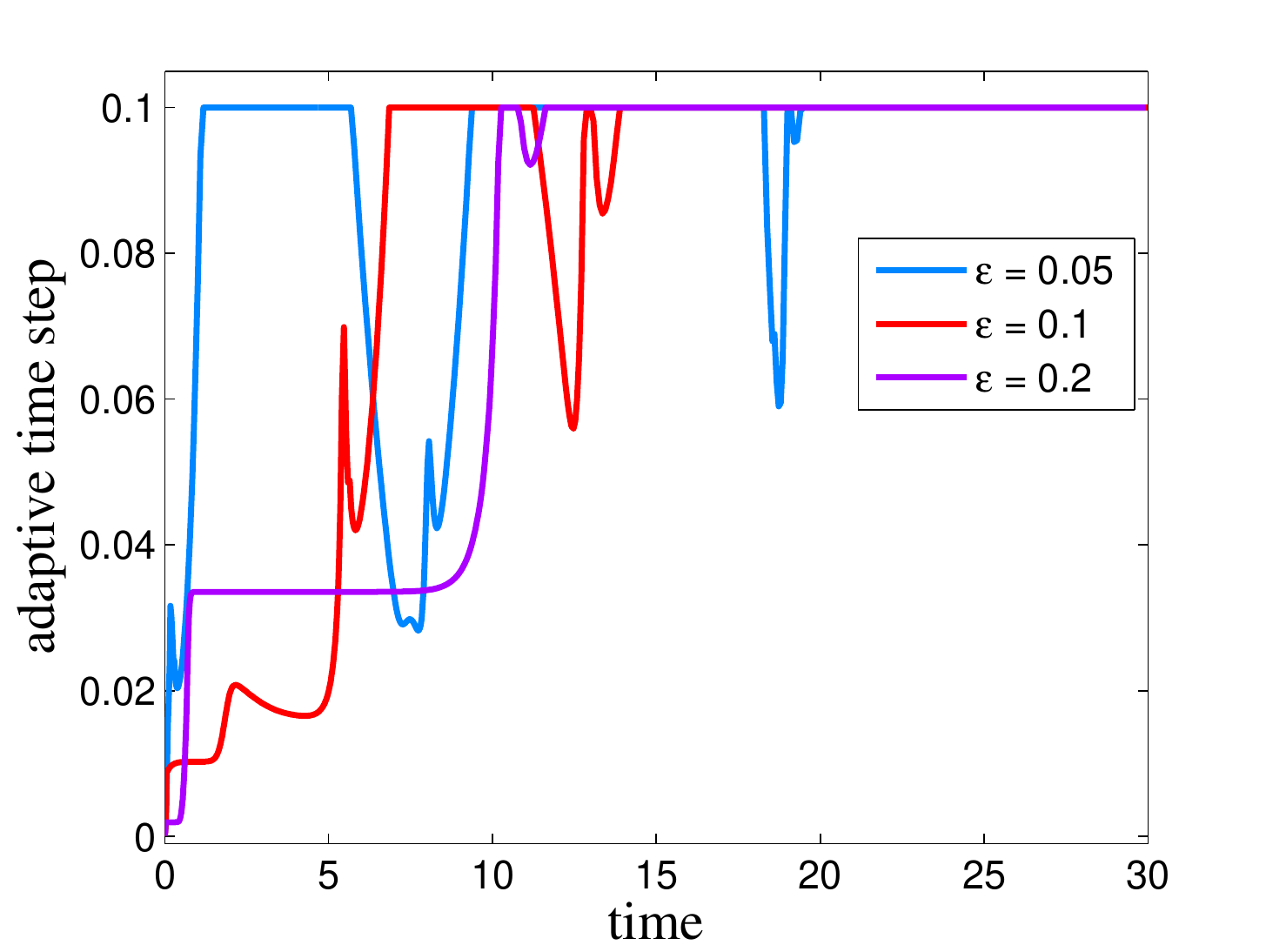}\\
\caption{Evolutions of energy (left) and
  time steps (right)  of  the MBE equation using initial data \eqref{example-2}
  with different $\varepsilon=0.2,0.1,0.05$
   until time T=30.}
\label{Comparison-epison-Energy}
\end{figure}

The evolutions of the phase variable obtained by adaptive time stepping strategy
are depicted in Figure \ref{example2-Snapshots-Dynamics}
and the evolution of the energy for the MBE model is presented
in Figure \ref{example2-Snapshots-Dynamics}.
The discrete energy, roughness, and adaptive time steps are shown in Figure \ref{example2-Adaptive-Energy-Roughness}.
In order to see the numerical performance, we use the same initial data with
different parameters $\epsilon=0.2,0.1,0.05$ to
carry out the simulations.
The energy curves and the correspondingly adaptive steps are summarized in
Figure \ref{Comparison-epison-Energy}.
We observe that the variable-step BDF2 scheme \eqref{scheme: BDF2 implicit MBE-no-slope}
 with the adaptive settings $\tau_{\max}=0.1$ and $\tau_{\min}=10^{-4}$
can work well for the current simulations.

\section{Conclusions}
We have performed the stability and convergence analysis of the variable-step BDF2 scheme for
the molecular beam epitaxial model without slope selection. The main contribution is that
we show that the variable-step BDF2 scheme admits an energy dissipation law
under the time-step ratios constraint $r_k:=\tau_k/\tau_{k-1}<3.561.$
Moreover, the $L^2$ norm stability and rigorous error estimates are established
under the same step-ratios constraint that ensuring the energy stability., i.e., $0<r_k<3.561.$
This is known to be the best result in literature. We remark that the technique
in this work is not applicable to molecular beam epitaxial model with slope selection,
and we shall pursuit this study in our future works.


\begin{thebibliography}{99}

\bibitem{AmarFamily:1996}
{\sc J.G. Amar and F. Family},
{\em Effects of crystalline microstructure on epitaxial growth},
 J Phys. Rev., 54 (1996), pp. 14071-14076.

\bibitem{Becker:1998}
{\sc J. Becker},
{\em A second order backward difference method with variable steps for a parabolic
problem}, BIT, 38(4) (1998), pp. 644--662.



 \bibitem{ChenCondeWangWangWise:2012}
{\sc W. Chen, S. Conde, C. Wang, X. Wang and S. M. Wise},
{\em A linear energy stable scheme for a thin film model without slope selection},
J Sci. Comput., 52 (2012), pp. 546-562.


\bibitem{ChenWangWang:2014}
{\sc W. Chen, C. Wang and X. Wang},
{\em A linear iteration algorithm for a second-order energy stable scheme for a thin film model without slope selection},
J Sci. Comput., 59 (2014), pp. 574--601.


\bibitem{ChenWangYanZhang:2019}
{\sc W. Chen, X. Wang, Y. Yan and Z. Zhang},
{\em A second order BDF numerical scheme with variable steps for the Cahn--Hilliard equation},
SIAM J. Numer. Anal., 57 (1) (2019), pp. 495--525.






\bibitem{CrouzeixLisbona:1984}
{\sc M. Crouzeix and F.J. Lisbona},
{\em The convergence of variable-stepsize, variable formula, multistep methods},
SIAM J. Numer. Anal., 21 (1984), pp. 512--534.





\bibitem{Emmrich:2005}
{\sc E. Emmrich},
{\em Stability and error of the variable two-step BDF for semilinear
parabolic problems}, J. Appl. Math. \& Computing, 19 (2005), pp. 33--55.



 \bibitem{EvansThiel:2010}
{\sc J.W. Evans, P.A. Thiel},
{\em A little chemistry helps the big get bigger},
Science, 330 (2010), pp. 599-600.





\bibitem{Grigorieff:1983}
{\sc R.D. Grigorieff},
{\em Stability of multistep-methods on variable grids}, Numer. Math., 42 (1983), pp. 359--377.


\bibitem{Golubovic:1997}
{\sc L. Golubovic},
{\em Interfacial coarsening in epitaxial growth models without slope selection},
Phys. Rev. Lett., 78 (1997), pp. 90-93.


\bibitem{GomezHughes:2011Provably}
{\sc H.~Gomez and T.~Hughes},
{\em Provably unconditionally stable, second-order time-accurate, mixed
  variational methods for phase-field models},
J. Comput. Phys., 230 (2011), pp.5310--5327.











\bibitem{JuLiQiaoZhang:2018}
{\sc L. Ju, X. Li, Z. Qiao and H. Zhang},
{\em Energy stability and error estimates of exponential time differencing schemes for the epitaxial growth model without slope selection},
Math. Comp., 87 (2018), pp. 1859--1885.



\bibitem{LeRoux:1982}
{\sc M.-N. Le Roux},
{\em Variable step size multistep methods for parabolic problems}, SIAM J.
Numer. Anal., 19 (4) (1982), pp. 725--741.
%







 \bibitem{LiLiu:2003}
{\sc B. Li and J.G. Liu},
{\em Thin film epitaxy with or without slope selection},
European J. Appl. Math., 14 (2003), pp. 713--743.

%



\bibitem{LiaoTangZhou:2019maximumAC}
{\sc H.-L. Liao, T. Tang and T. Zhou},
{\em On energy stable, maximum-principle preserving, second order
BDF scheme with variable steps for the Allen-Cahn equation},
SIAM J. Numer. Anal., 2020, to appear.

 \bibitem{LiaoZhang:2019linear}
{\sc H.-L. Liao and Z. Zhang},
{\em Analysis of adaptive BDF2 scheme for diffusion equations},
Math. Comp., 2020, to appear.


\bibitem{LiaoJiZhang:2020}
{\sc H.-L Liao,  B. Ji and L. Zhang},
{\em An adaptive BDF2 implicit time-stepping method for the phase field crystal model},
IMA J. Numer. Anal., 2020, in review (arXiv:2008.00212v1).










\bibitem{QiaoSunZhang:2015}
{\sc Z. Qiao, Z.Z. Sun and Z. Zhang},
{\em Stability and convergence of second-order schemes for the nonlinear epitaxial growth model without slope selection},
Math Comp., 84 (2015), pp. 653-674.



\bibitem{QiaoZhangTang:2011}
{\sc Z. Qiao, Z. Zhang and T. Tang},
{\em An adaptive time-stepping strategy for the molecular beam epitaxy models},
SIAM J. Sci. Comput., 33 (2011), pp. 1395--1414.
%





\bibitem{RostKrug:1997}
{\sc M. Rost and J. Krug},
{\em Coarsening of surfaces in unstable epitaxial growth},
J Phys. Rev. E., 55 (1997), pp. 4952-3957.





 \bibitem{ShenWangWangWise:2012}
{\sc J. Shen, C. Wang, X. Wang and S.M. Wise},
{\em Second-order convex splitting schemes for gradient flows with
Ehrlich--Schwoebel type energy: application to thin film epitaxy},
SIAM J. Numer. Anal., 50(1) (2012), pp. 105--125.



\bibitem{XuLiWuBousquet:2019}
{\sc J. Xu, Y.K. Li, S.N. Wu and A.Bousequet},
{\em On the stability and accuracy of partially and fully implicit schemes for phase field modeling},
Comput. Methods Appl. Mech. Engrg., 345 (2019), pp. 826-853.


\bibitem{XuTang:2006}
{\sc C. Xu and T. Tang},
{\em Stability analysis of large time-stepping methods for epitaxial growth models},
SIAM J. Numer. Anal. 44(4) (2006), pp. 1759--1779.







\bibitem{ZhangMaQiao:2013}
{\sc Z. Zhang,  Y. Ma and Z. Qiao},
{\em An adaptive time-stepping strategy for solving the phase field crystal model},
 J. Comput. Phys., 249 (2013), pp. 204--215.










\end{thebibliography}
\end{document}